\documentclass{siamltex}
\usepackage{amsfonts}
\usepackage{amssymb,amsmath}
\usepackage{graphicx}
\usepackage{subfigure}
\usepackage{color}
\usepackage{fullpage}

\DeclareMathAlphabet{\itbf}{OML}{cmm}{b}{it}

\def\EE{\mathbb{E}}

\def\RR{\mathbb{R}}
\def\eps{\varepsilon}

\def\om{{\omega}}
\def\la{{\lambda}}

\def\cA{{\mathcal A}}
\def\cB{{\mathcal B}}

 \def\cS{{\mathcal S}}
\def\cT{{\mathcal T}}

\def\cO{{\mathcal O}}
\def\cD{{\mathcal D}}

\def\txi{\widetilde{\xi}}

\def\tk{\widetilde{\kappa}}
\def\tom{\widetilde{\omega}}
\renewcommand{\hat}{\widehat}
\newcommand{\IKM}{{\cal J}^{\tiny \mbox{M}}}
\newcommand{\ICINT}{{\cal J}^{\tiny \mbox{CINT}}}
\newcommand{\WICINT}{\widetilde{\cal J}^{\tiny \mbox{CINT}}}
\newcommand{\ITR}{{\cal J}^{\tiny \mbox{TR}}}
\newcommand{\DTR}{D^{\tiny \mbox{TR}}}
\newcommand{\hDTR}{\widehat{D}^{\tiny \mbox{TR}}}
\newcommand{\FTR}{{\mathcal M}^{\tiny \mbox{TR}}}
\newcommand{\FKM}{{\mathcal M}^{\tiny \mbox{M}}}
\newcommand{\FCINT}{{\mathcal M}^{\tiny \mbox{CINT}}}
\newcommand{\WFCINT}{\widetilde{\mathcal M}^{\tiny \mbox{CINT}}}

\begin{document}

\title{Paraxial coupling of propagating modes in three-dimensional
  waveguides with random boundaries}

\author{  Liliana Borcea\footnotemark[1] and Josselin
Garnier\footnotemark[2] }

\maketitle

\renewcommand{\thefootnote}{\fnsymbol{footnote}}

\footnotetext[1]{Computational    and   Applied    Mathematics,   Rice
  University,    Houston,   TX    77005.   {\tt    
  borcea@rice.edu}} \footnotetext[2]{Laboratoire  de Probabilit\'es et
  Mod\`eles   Al\'eatoires   \&   Laboratoire   Jacques-Louis   Lions,
  Universit{\'e}  Paris VII,  Site Chevaleret,  75205 Paris  Cedex 13,
  France.  {\tt garnier@math.jussieu.fr}}

\renewcommand{\thefootnote}{\arabic{footnote}}

\begin{abstract}
  We analyze long range wave propagation in three-dimensional random
  waveguides. The waves are trapped by top and bottom boundaries, but
  the medium is unbounded in the two remaining directions. We consider
  scalar waves, and motivated by applications in underwater acoustics,
  we take a pressure release boundary condition at the top surface and
  a rigid bottom boundary. The wave speed in the waveguide is known
  and smooth, but the top boundary has small random fluctuations that
  cause significant cumulative scattering of the waves over long
  distances of propagation.  To quantify the scattering effects, we
  study the evolution of the random amplitudes of the waveguide modes.
  We obtain that in the long range limit they satisfy a system of
  paraxial equations driven by a Brownian field.  We use this
  system to estimate three important mode-dependent scales: the
  scattering mean free path, the cross-range decoherence length and
  the decoherence frequency. Understanding these scales is important
  in imaging and communication problems, because they encode
  the cumulative scattering effects in the wave field measured by
  remote sensors.  As an application of the theory, we analyze time
  reversal and coherent interferometric imaging in strong cumulative
  scattering regimes.
\end{abstract}



\section{Introduction}
\label{sect:intro}
\begin{figure}[t!]
\vspace{-0.4in}
   \begin{center}
     \input{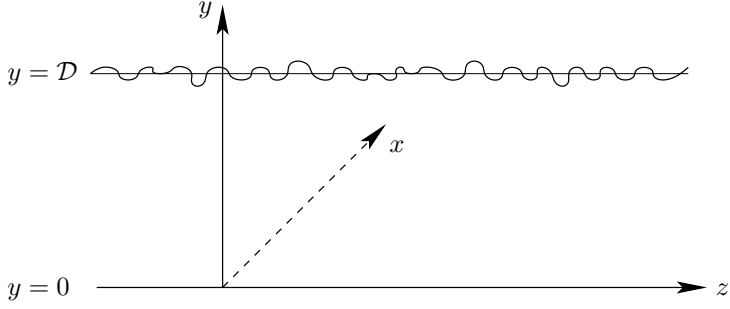}
   \end{center}
   \vspace{-0.05in}
   \caption{Schematic of the problem setup. The system of coordinates
    has range origin $z=0$ at the source. The rigid bottom boundary
    $y=0$ is assumed flat and the pressure release top boundary has 
    fluctuations around the value $y = {{\mathcal D}}$. The cross-range $x$ and the 
    range $z$ are unbounded, that is $(x,z) \in \mathbb{R}^2$.} 
   \label{fig:schematic}
\end{figure}

We study long range scalar (acoustic) wave propagation
 in a three-dimensional waveguide.  The setup is illustrated in Figure
\ref{fig:schematic}, and it is motivated by problems in underwater
acoustics. We denote by $z \in \mathbb{R}$ the range, the main
direction of propagation of the waves.  The medium is unbounded in the
cross-range direction $x \in \mathbb{R}$, but it is confined in depth
$y$ by two boundaries which trap the waves, thus creating the
waveguide effect.

The acoustic pressure field is denoted by $p(t,x,y,z)$, and it
satisfies the wave equation 
\begin{equation}
  \left[ \partial_x^2 + \partial_y^2 + \partial_z^2 - 
   \frac{1}{c^2(y)} \partial_t^2 \right] p(t,x,y,z) = f(t,x,y)\delta(z), 
\qquad y \in [0,T(x,z)], ~ x, z \in \mathbb{R}, ~ ~ t > 0,
\label{eq:form.1}
\end{equation}
in a medium with wave speed $c(y)$. The excitation is due to a source
located in the plane $z=0$, emitting the pulse $f(t)$.  The medium is
quiescent before the source excitation,
\begin{equation}
p(t,x,y,z) = 0, \quad t \ll 0.
\label{eq:form.2}
\end{equation}
The bottom of the waveguide is assumed rigid
\begin{equation}
\partial_y p(t,x,y=0,z) = 0,
\label{eq:form.3}
\end{equation}
and we take a pressure release boundary condition at the perturbed top
boundary
\begin{equation}
p(t,x,y=T(x,z),z) = 0.
\label{eq:form.4}
\end{equation}
Perturbed means that the boundary $y = T(x,z)$ has small fluctuations
around the mean depth ${{\mathcal D}}$,
\begin{equation}
|T(x,z) - {{\mathcal D}}| \ll {{\mathcal D}}.
\label{eq:form.5}
\end{equation}
We choose this setup for simplicity. The results extend readily to
other boundary conditions and to fluctuating bottoms.  Such
boundaries were considered recently in \cite{ABG-12,gomez2}, 
in two-dimensional waveguides.  Extensions to media with small 
$(x,y,z)$-dependent fluctuations of the wave speed can also be made using the
techniques developed in
\cite{kohler77,dozier,book07,garnier_papa,gomez}.

The goal of our study is to quantify the effect of scattering at the
surface. Because in applications it is not feasible to know the
boundary fluctuations in detail, we model them with a random process.
The solution $p(t,x,y,z)$ of equations
(\ref{eq:form.1})-(\ref{eq:form.4}) is therefore a random field, and
we describe in detail its statistics at long ranges, where cumulative
scattering is significant. We use the results for two applications:
time reversal and sensor array imaging.

Our method of solution uses a change of coordinates to straighten the
boundary. The transformed problem has a simple geometry but a randomly
perturbed differential operator.  Its solution is given by a
superposition of propagating and evanescent waveguide modes, with
random amplitudes. We show that in the long range limit these
amplitudes satisfy a system of paraxial equations that are driven by a
Brownian field.  The detailed characterization of the statistics of
$p(t,x,y,z)$ follows from this system. It involves the calculation of
the mode-dependent scattering mean free path, which is the distance
over which the modes lose coherence; the mode-dependent decoherence
length, which is the cross-range offset over which the mode amplitudes
decorrelate; and the mode-dependent decoherence frequency, which is
the frequency offset over which the mode amplitudes decorrelate. These
scales are important in studies of time reversal and imaging, because
they dictate the resolution of focusing and the robustness
(statistical stability) of the results with respect to realizations of
the random fluctuations of the boundary.

The paper is organized as follows: We begin in section
\ref{sect:HOMWG} with the description of the reference pressure field
$p_o(t,x,y,z)$ in ideal waveguides with planar boundaries.  The random
field $p(t,x,y,z)$ derived in section \ref{sect:RANDWG} may be viewed
as a perturbation of $p_o(t,x,y,z)$, in the sense that it is
decomposed in the same waveguide modes. However, the amplitudes of the
modes are random and coupled. Because the fluctuations of the boundary
are small, we consider in section \ref{sect:LR} a long range limit, so
that we can observe significant cumulative scattering.  The statistics
of the wave field at such long ranges is described in section
\ref{sect:STAT}. The results are summarized in section \ref{sect:FM}
and are used in sections \ref{sect:TR} and \ref{sect:IM} for analyzing
time reversal and imaging with sensor arrays. We end with a summary in
section \ref{sect:SUM}.

\section{Wave propagation in ideal waveguides}
\label{sect:HOMWG}
The pressure field in ideal waveguides, with planar boundaries, is
given by
\begin{equation}
  p_o(t,x,y,z) = \int_{-\infty}^\infty \frac{d \om}{2 \pi} \, 
\hat p_o(\om,x,y,z)  e^{- i \om t},
\label{eq:HOM.1}
\end{equation}
with Fourier coefficients satisfying a separable problem for the
Helmholtz equation
\begin{equation}
  \left[ \partial_x^2 + \partial_y^2 + \partial_z^2 + \frac{\om^2}{c^2(y)} 
  \right] \hat p_o(\om,x,y,z) = \widehat f(\om,x,y) \delta(z), \qquad 
  |\om-\om_0| \le \frac{B}{2}, ~ ~ (x,z) \in \mathbb{R}^2, ~~ 
  y \in (0,{{\mathcal D}}),
\label{eq:HOM.2}
\end{equation}
with boundary conditions 
\begin{equation}
\partial_y \hat p_o(\om,x,y=0,z)= \hat p_o(\om,x,y={{\mathcal D}},z) = 0,
\label{eq:HOM.3}
\end{equation}
and outgoing radiation conditions at $\sqrt{x^2 + z^2} \to \infty$.
The Fourier transform of the source
\begin{equation}
\hat f(\om,x,y) = \int_{-\infty}^\infty d t \, f(t,x,y) e^{i \om t},
\end{equation}
is compactly supported in $[\om_0- B/2,\om_0+  B/2]$, for any $x$ and
$y$. Here $\om_0$ is the central frequency and $B$ is the bandwidth.

\subsection{Propagating and evanescent modes}
\label{sect:form.Prop}
The solution of the Helmholtz equation (\ref{eq:HOM.2}) is a
superposition of $N(\om)$ propagating modes, and infinitely many
evanescent ones,
\begin{equation}
  \hat p_o(\om,x,y,z) = \sum_{j=1}^{N(\om)} \phi_j(\om,y) 
\hat u_{j,o}(\om,x,z)  
  + \sum_{j = N(\om)+1}^{\infty} \phi_j(\om,y) \hat v_{j,o}(\om,x,z).
\label{eq:HOM.4}
\end{equation}
The decomposition is in the $L^2(0,{{\mathcal D}})$ orthonormal basis of
the eigenfunctions $\phi_j(\om,y)$ of the self-adjoint differential
operator in $y$,
\begin{eqnarray*}
  \left[ \partial_y^2 + \frac{\om^2}{c^2(y)}\right] \phi_j(\om,y) &=&
  \lambda_j(\om) \phi_j(\om,y), \nonumber \\ \phi_j(\om,{{\mathcal
  D}}) &=& \partial_y \phi_j(\om,0) = 0, \qquad j = 1, 2, \ldots,
\end{eqnarray*}
with eigenvalues $\lambda_j(\om)$ that are simple \cite{Weid}.

To simplify the analysis, we assume in this paper that the wave speed
is homogeneous
\begin{equation}
c(y) = c_o.
\label{eq:HOM.6}
\end{equation}
The results for variable wave speeds are similar in all the essential
aspects.  The simplification brought by (\ref{eq:HOM.6}) amounts to
having explicit expressions of the eigenfunctions, which are
independent of the frequency
\begin{equation}
\phi_j(y) = \sqrt{\frac{2}{{{\mathcal D}}}} \cos \left[ \pi
\Big(j-\frac{1}{2}\Big) \frac{y}{{{\mathcal D}}} \right].
\label{eq:HOM.7}
\end{equation}
The eigenvalues are
\begin{equation}
  \lambda_j(\om)= \left(\frac{\pi}{{{\mathcal D}}}\right)^2 \left[ 
\Big(\frac{k {{\mathcal D}}}{\pi}\Big)^2 - \Big(j-\frac{1}{2}\Big)^2\right],
\label{eq:HOM.8}
\end{equation}
where $k = \om/c_o$ is the wavenumber, and only the first $N(\om)$ of
them are non-negative
\begin{equation}
N(\om) = \left \lfloor \frac{k {{\mathcal D}}}{\pi} + \frac{1}{2}
\right \rfloor.
\label{eq:HOM.9}
\end{equation}
The notation $\lfloor ~ \rfloor$ stands for the integer part.
We suppose for simplicity that $N(\om)$ remains constant in the 
bandwidth $[\omega_0-B/2,\omega_0+B/2]$, and write from now on $N(\om) = N$. 

The propagating components in (\ref{eq:HOM.4}) satisfy 
the two-dimensional Helmholtz equation
\begin{equation}
  \left[ \partial_x^2 + \partial_z^2 + \beta_j^2(\om) \right] 
\hat u_{j,o}(\om,x,z) = 
  \hat F_j(\om,x) \delta(z), \qquad j = 1, \ldots, N,
\label{eq:HOM.10}
\end{equation}
with outgoing, radiation conditions at $\sqrt{x^2 + z^2} \to \infty$.
The evanescent components solve
\begin{equation}
  \left[ \partial_x^2 + \partial_z^2 - \beta_j^2(\om) \right] 
  \hat v_{j,o} (\om,x,z)=   \hat F_j(\om,x) \delta(z), \qquad j > N,
\label{eq:HOM.11}
\end{equation}
with decay condition $ \hat v_{j,o}(\om,x,z) \to 0$ at $\sqrt{x^2 +
  z^2} \to \infty.$ Here we introduced the coefficients of the source profile in the basis of the
eigenfunctions
\begin{equation}
  \hat F_j(\om,x) = \int_0^{{\mathcal D}} dy \, \phi_j(y) \hat f(\om,x,y), 
  \qquad j \ge 1, 
\label{eq:HOM.14}
\end{equation}
and the mode wavenumbers
\begin{equation}
  \beta_j(\om) = \sqrt{|\lambda_j(\om)|} = \frac{\pi}{{{\mathcal D}}}
  \sqrt{ \left|\Big(\frac{k {{\mathcal D}}}{\pi}\Big)^2 -
  \Big(j-\frac{1}{2}\Big)^2\right|}, \qquad j \ge 1.
\label{eq:HOM.13}
\end{equation}
We assume that none of the $\beta_j(\om)$  vanishes in the 
bandwidth, so that there are no standing waves. That is to say, 
\begin{equation}
  \frac{k{{\mathcal D}}}{\pi} = N + \alpha(\om) -\frac{1}{2}, \qquad 
  \alpha(\om) \in (0,1)
  ~~  \mbox{for all} ~~ \om \in [\om_0- B/2,\om_0+ B/2].
\label{eq:DefAlpha}
\end{equation}

\subsection{The paraxial regime}
We now introduce the paraxial scaling for the ideal waveguide, with
the source emitting a beam that propagates along the $z$-axis.  As we
show below, this happens when the cross-range profile of the source is
larger than the wavelength.  The source generates a quasi-plane wave,
with slowly varying envelope satisfying a Schr\"odinger-like equation.

Explicitly, we assume that the source is of the form
\begin{equation}
\label{scaledsource}
f^\eps(t,x,y) = f(t,\eps x ,y)
\end{equation}
where $\eps$ is a small dimensionless parameter defined as the ratio
of the central wavelength $\lambda_{0}$ and the transverse width
$r_{0}$ of the source.  Standard diffraction theory gives that the
Rayleigh length for a beam with initial width $r_0 = \lambda_{0}/\eps$
is of the order of
\[
r_0^2 /\lambda_0 = \lambda_{0}/\eps^{2}.
\] 
The Rayleigh length is the distance along the $z$ axis from the beam
waist to the place where the beam area is doubled by diffraction.
Therefore, we look at the wavefield at $ O(\eps^{-1})$ cross-range
scales, similar to $r_0$, and at $O(\eps^{-2})$ range scale, similar
to the Rayleigh length. We rename the field in this scaling as
 \begin{equation}
 \label{eq:poeps}
p_o^\eps (t,X,y,Z) = p_o \Big(t
,\frac{X}{\eps},y,\frac{Z}{\eps^2}\Big).
\end{equation}

The Fourier coefficients of (\ref{eq:poeps}) are given by the scaled
version of (\ref{eq:HOM.4})
\begin{equation}
  \hat p_o^\eps(\om,X,y,Z) = \sum_{j=1}^{N(\om)} \phi_j (y) 
\hat u_{j,o}^\eps(\om,X,Z)  
  + \sum_{j = N(\om)+1}^{\infty} \phi_j(y) \hat v_{j,o}^\eps(\om,X,Z),
\label{eq:HOM.4eps}
\end{equation}
with propagating mode amplitudes $\hat u_{j,o}^\eps$ satisfying the
scaled equation (\ref{eq:HOM.10}), with the source replaced by $\hat
F_{j}(\om,\eps x = X)$. They can be written as
\[
\hat u_{j,o}^\eps(\om,X,Z) =- \frac{1}{\eps} \int_{-\infty}^\infty d
X' \, \hat F_j(\om,X') \hat G_o \Big(\beta_j(\om),
\frac{X-X'}{\eps},\frac{Z}{\eps^2}\Big),
\]
in terms of the outgoing Green's function 
\[
\hat G_o \big(\beta_j(\om),x,z \big) = \frac{i}{4}
H_0^{(1)}\left[\beta_j(\om) \sqrt{x^2 + z^2}\right].
\]
Here $H_0^{(1)}$ is the Hankel function of the first kind, and because
$\eps \ll 1$, we can use its asymptotic form for a scaled range $Z>0$
\begin{eqnarray*}
  \frac{i}{4}H_0^{(1)} \left[ \beta_j(\om)
  \sqrt{\frac{(X-X')^2}{\eps^2} + \frac{Z^2}{\eps^4}}\right] &\approx&
  \frac{1}{4}\left[ \frac{2 i}{\pi \beta_j(\om) \sqrt{
  \frac{(X-X')^2}{\eps^2} + \frac{Z^2}{\eps^4}}} \right]^{1/2} \exp  \left[i
  \beta_j(\om) \sqrt{ \frac{(X-X')^2}{\eps^2} + \frac{Z^2}{\eps^4}}\right]
  \\ &\approx& \frac{\eps}{2 }\sqrt{ \frac{i}{2 \pi \beta_j(\om) Z} }
  \exp \left\{i \beta_j(\om) \left[\frac{Z}{\eps^2} + \frac{(X-X')^2}{2
  Z}\right]\right\}.
\end{eqnarray*}
The propagating components of the wave field become
\[
\hat u_{j,o}^\eps(\om,X,Z) \approx a_{j,o}(\om,X,Z) \exp\left[i
\beta_j(\om) \frac{Z}{\eps^2}\right],
\]
with 
\begin{equation}
  a_{j,o}(\om,X,Z) = -\frac{1}{2}\sqrt{ \frac{i}{2 \pi \beta_j(\om) Z
  } } \int_{-\infty}^\infty d X' \, \exp\left[ \frac{i \beta_j(\om)
  (X-X')^2}{2 Z}\right] \hat F_j(\om,X'),
\label{eq:HOM.15}
\end{equation}
for $j = 1, \ldots, N$. The evanescent components are obtained
similarly from (\ref{eq:HOM.11})
\[
\hat v_{j,o}^\eps(\om,X,Z) \approx e_{j,o}(\om,X,Z) \exp
\left[-\beta_j(\om) \frac{Z}{\eps^2}\right],\] and
\begin{equation}
e_{j,o}(\om,X,Z) = -\frac{1}{2}\sqrt{ \frac{1}{2 \pi \beta_j(\om) Z }
  } \int_{-\infty}^\infty d X' \, \exp\left[ -\frac{ \beta_j(\om)
  (X-X')^2}{2 Z}\right] \hat F_j(\om,X'),
\label{eq:HOM.17}
\end{equation}
for $j \ge N+ 1$. These modes are exponentially damped and can be neglected.

In summary, the paraxial approximation of the wave field is given by 
\begin{equation}
  \hat p_o^\eps(\om,X,y,Z) \approx \sum_{j=1}^{N} \phi_j(y)
a_{j,o}(\om,X,Z) e^{i \beta_j(\om) \frac{Z}{\eps^2}}
\label{eq:HOM.parax}
\end{equation}
It is a superposition of forward going modes with  complex
valued amplitudes $a_{j,o}$
given by (\ref{eq:HOM.15}), and solving the paraxial equations
\begin{equation}
  \left[ 2 i \beta_j(\om) \partial_Z + \partial_X^2\right] a_{j,o}(\om,X,Z) 
= 0 , \qquad j = 1, \ldots, N, 
\label{eq:HOM.16}
\end{equation}
with initial conditions
\begin{equation}
a_{j,o}(\om,X,Z=0) = a_{j,{\rm ini}} (\om,X):= \frac{1}{2 i
\beta_j(\omega)} \hat F_j(\om,X) , \qquad j = 1, \ldots, N.
\label{eq:HOM.16b}
\end{equation}

\section{Wave propagation in random waveguides}
\label{sect:RANDWG}
In this section we consider a waveguide with fluctuating boundary and
analyze the wave field under the following scaling assumptions:
\begin{enumerate}
\item The
transverse width $r_0$ of the source and the central wavelength
$\lambda_0$ satisfy
\begin{equation}
\label{eq:sc.1}
r_0 = \eps^{-1} \lambda_0,
\end{equation}
as in the previous section.  The correlation length $\ell^\eps$ of the
boundary fluctuations is similar to $r_0$,
\begin{equation}
\ell^\eps = \eps^{-1} \ell \sim r_0,
\label{eq:sc.2}
\end{equation}
so that there is a non-trivial interaction between the boundary
fluctuations and the wavefield. Here $\ell$ is the scaled order-one
correlation length defined below.
\item  The scale $L^\eps$ of the
propagation distance is much larger than $\lambda_0$.  More precisely,
\begin{equation}
\label{eq:sc.3}
L^\eps/\lambda_0 = O(\eps^{-2}).
\end{equation}
Recall that the Rayleigh length for a beam with initial width $r_0$
and central wavelength $\lambda_0$ is of the order of $r_0^2
/\lambda_0 \sim \eps^{-2} \lambda_0$ in absence of random
fluctuations.  The high-frequency scaling assumption (\ref{eq:sc.3})
ensures that the propagation distance is similar to the Rayleigh
length.
\item  The amplitude of the boundary
fluctuations is small, of the order of $\eps^{3/2} \lambda_0$.  As we
will show, this scaling is precisely the one that gives a cumulative
scattering effect of order one after the propagation distance
$L^\eps$.
\end{enumerate}

\vspace{0.05in} We use the hyperbolicity of the problem to truncate
mathematically the boundary fluctuations to the range interval
$(0,L/\eps^2)$. The bound $L/\eps^2$ is the maximum range of the
fluctuations that can affect the waves up to the observation time
$T^\eps$ of order $\eps^{-2}$. The lower bound in the range
interval coincides with the location of the source. It is motivated by
two facts: First, we observe the waves at positive ranges. Second, the
backscattered field is negligible in the scaling regime defined above,
as we show later in section \ref{sect:LR.2}.

The boundary fluctuations are modeled with a random process $\mu$
\begin{equation}
  T^\eps(x,z) = {{\mathcal D}} \left[ 1 + \eps^{3/2} \mu\left(\eps x, 
\eps z \right)\right], 
    \qquad z \in (0,  L/\eps^2).
\label{eq:RAND.1}
\end{equation}
The process $\mu$ is bounded, zero-mean, stationary and mixing,
meaning in particular that its covariance is integrable\footnote{More
  precisely, $\mu$ is a $\varphi$-mixing process, with $\varphi \in
  L^{1/2}(\RR^+)$, as stated in \cite[4.6.2]{kushner}.}. Because our
method of solution flattens the boundary by changing coordinates, we
require that $\mu$ is twice differentiable, with almost surely bounded
derivatives.
 Its covariance function is given by
\begin{equation}
  R\left(\xi,\zeta\right) = 
  \EE\left[ \mu(\xi'+\xi,\zeta'+\zeta)\mu(\xi',\zeta')\right],
\label{eq:RAND.2}
\end{equation}
and we denote by $R_o(\xi)$ its integral over $\zeta$,
\begin{equation}
  R_o(\xi) = \int_{-\infty}^\infty d \zeta \, R(\xi,\zeta) .
\label{eq:RAND.3p}
\end{equation}
Our assumption on the differentiability of $\mu$ implies that $R_o$ is
four times differentiable.  Note that $\xi = 0$ is the maximum of the
integrated covariance $R_o(\xi)$, so we have
\begin{equation}
R_o'(0) = 0.
\label{eq:RAND.4}
\end{equation}
We define the scaled square amplitude $\sigma^2$ and correlation
length $\ell$ of the boundary fluctuations through the equations
\begin{equation}
R_o(0) = \sigma^2 \ell, 
 \quad \quad \frac{R_o''(0)}{R_o(0)} = -\frac{1}{\ell^2} .
\label{eq:RAND.3}
\end{equation}

\subsection{Change of coordinates}
\label{sect:RANDWG.2}
We introduce the change of coordinates from $(x,y,z)$ to $(x,\eta,
z)$, with
\begin{equation}
\eta = \frac{y {{\mathcal D}}}{T^\eps(x,z) }.
\end{equation}
It straightens the boundary $y = T^\eps (x,z)$ to $\eta = {{\mathcal
D}}$, for any $x \in \mathbb{R}$ and $z \in (0,L/\eps^2)$. The
pressure field in the new coordinates is denoted by
\begin{equation}
\hat P (\om,x,\eta,z) = \hat p\left(\om,x,\frac{\eta
 T^\eps(x,z)}{{{\mathcal D}}}, z \right).
\label{eq:RAND.6}
\end{equation}
It satisfies the simple  boundary conditions  
\begin{equation}
\hat P (\om,x,{{\mathcal D}},z) = \partial_\eta \hat P (\om,x,0,z) = 0,
\label{eq:RAND.7}
\end{equation}
and the partial differential equation 
\begin{eqnarray}
\left[ \partial_x^2 + \partial_z^2 + \left(\frac{{{\mathcal
D}}^2}{{T^\eps }^2} + \eta^2 \frac{|\nabla T^\eps |^2}{{T^\eps
}^2}\right) \partial_\eta^2 - 2 \eta \frac{\nabla T^\eps }{T^\eps }
\cdot \nabla \partial_\eta + \left( 2 \eta \frac{ |\nabla T^\eps
|^2}{{T^\eps }^2} - \eta \frac{\Delta T^\eps }{T^\eps }\right)
\partial_\eta + k^2 \right] \hat P = \nonumber \\
\hat{f}^\eps (\omega,x, \eta )
\delta(z),
\label{eq:RAND.8}
\end{eqnarray}
derived from (\ref{eq:form.1}) and (\ref{eq:RAND.6}) using the chain rule.
Here $\nabla$ and $\Delta$ are the gradient and Laplacian operators in
$(x,z)$ and $f^\eps$ is the source of the form (\ref{scaledsource}).  

When substituting the model (\ref{eq:RAND.1}) in (\ref{eq:RAND.8}), we
obtain that $\hat P$ satisfies a randomly perturbed problem
\begin{equation}
  \left[ \partial_x^2 + \partial_z^2 +\left(1-2 \eps^{3/2} \mu(\eps
x,\eps z) \right) \partial_\eta^2 + k^2 + \ldots \right] \hat
P(\om,x,\eta,z) = \hat{f} (\omega, \eps x, \eta ) \delta(z) .
\label{eq:RAND.9}
\end{equation}
The higher-order terms denoted by the dots are
\[
-2 \eps^{5/2} \left[1 + O(\eps^{3/2}) \right] \eta \nabla \mu \cdot
\nabla \partial_\eta \hat P +
3 \eps^3 \mu^2 \left[1+O(\eps^2)\right] \partial_\eta^2 \hat P -
\eps^{7/2} \Delta \mu \left[1 + O(\eps^{3/2})\right]\eta
 \partial_\eta \hat P.
\]
They come from the  expansions in $\eps$ of the coefficients 
in (\ref{eq:RAND.8}), and are negligible in the limit $\eps \to 0$
considered in section \ref{sect:LR}.

\subsection{Wave decomposition}
\label{sect:RANDWG.3}
Equation (\ref{eq:RAND.9}) is not separable, but we can still write
its solution in the $L^2(0,{{\mathcal D}})$ basis of the eigenfunctions
(\ref{eq:HOM.7}). The
expansion is similar to (\ref{eq:HOM.4})
\begin{equation}
  \hat P(\om,x,\eta,z) = \sum_{j=1}^N \phi_j(\eta) 
\hat u_j(\om,x,z)  + \sum_{j>N} \phi_j(\eta) \hat v_j(\om,x,z).
\label{eq:RAND.10}
\end{equation}
We define the forward and backward going wave mode amplitudes $a_j$
and $b_j$ by
\begin{eqnarray}
a_j(\omega,x,z) &=& \Big( \frac{1}{2} \hat u_j(\om,x,z) + \frac{1}{2 i
\beta_j(\om)} \partial_z \hat u_j(\om,x,z)\Big) e^{- i \beta_j(\om) z}
, \nonumber \\ b_j(\omega,x,z) &=& \Big( \frac{1}{2} \hat u_j(\om,x,z) -
\frac{1}{2 i \beta_j(\om)} \partial_z \hat u_j(\om,x,z) \Big) e^{ i
\beta_j(\om) z}, \label{eq:def.ab}
\end{eqnarray}
so that the complex valued amplitudes of the propagating modes can be
written as
$$
\hat u_j(\om,x,z)  = a_j(\om,x,z) e^{i \beta_j(\om) z} + b_j(\om,x,z) 
e^{-i \beta_j(\om) z}.
$$ 

Definition (\ref{eq:def.ab}) implies that
\begin{equation}
\partial_z a_j(\om,x,z) e^{i \beta_j(\om) z} + \partial_z 
b_j(\om,x,z) e^{-i \beta_j(\om) z} = 0, \qquad j = 1, \ldots, N.
\label{eq:RAND.13}
\end{equation}
This equation is needed to specify uniquely the propagating mode
amplitudes, because they each satisfy a single boundary condition in
the range $(0,L/\eps^2)$ of the fluctuations.  To derive these boundary
conditions, let us observe that $a_j$ and $b_j$ must be constant in
$z\in (-\infty,0)$ and in $z \in (L/\eps^2,\infty)$, because the
boundary is flat outside $(0,L/\eps^2)$.  Moreover, the radiation
conditions
\[
\lim_{z \to -\infty} a_j(\om,x,z) = 0, \quad 
\lim_{z \to \infty} b_j(\om,x,z) = 0,
\]
imply that the mode amplitudes satisfy
\begin{eqnarray}
a_j(\omega,x,z=0^-) =0, \\ 
b_j(\omega,x,z=L/\eps^2) =0 .
\label{eq:RAND.12}
\end{eqnarray}
The last equation is the boundary condition for $b_j$.
The boundary value $a_j(\om,x, z=0^+)$ follows from the jump conditions 
across the plane $z = 0$ of the source in 
equation (\ref{eq:RAND.9}). We have 
\[
[ \hat{u}_j]_{0^-}^{0^+} =0, \quad [\partial_z \hat{u}_j]_{0^-}^{0^+} =\hat{F}_j(\omega,\eps x),
\]
with $\widehat{F}_j$ defined by (\ref{eq:HOM.14}).
This gives
\[[ \hat{a}_j + \hat{b}_j ]_{0^-}^{0^+} =0, \quad i \beta_j [ \hat{a}_j - \hat{b}_j]_{0^-}^{0^+} =\hat{F}_j(\omega,\eps x),\]
and therefore
\begin{equation}
  a_j(\om,x,0^+) = \frac{1}{2 i \beta_j(\omega)}
  \widehat{F}_j(\omega,\eps x) .
\label{eq:RAND.11}
\end{equation}

Substituting (\ref{eq:RAND.10}) in (\ref{eq:RAND.9}), and using the
orthonormality of the eigenfunctions $\phi_j$, we find that the wave
mode amplitudes solve paraxial equations coupled by the random
fluctuations in $z \in (0,L/\eps^2)$,
\begin{eqnarray}
  \left( 2 i \beta_j \partial_z + 
    \partial_x^2 \right) a_j + e^{-2 i \beta_j z} \partial_x^2 b_j &\approx& 
  \eps^{3/2}   \mu(\eps x, \eps z) e^{-i \beta_j z} \left[\sum_{l=1}^N q_{jl}
    \left( a_l e^{i \beta_l  z} + b_l 
      e^{-i \beta_l  z}\right) + 
     \sum_{l>N} q_{jl} \hat v_l \right], ~ ~ ~
  \label{eq:RAND.14} \\ 
  \left(- 2 i \beta_j \partial_z + 
    \partial_x^2 \right) b_j + e^{2 i \beta_j z}\partial_x^2 a_j & \approx &
  \eps^{3/2}   \mu(\eps x,\eps z)    e^{i \beta_j z}
  \left[\sum_{l=1}^N q_{jl}
    \left( a_l e^{i \beta_l z} + b_l 
      e^{-i \beta_l  z}\right) + \sum_{l>N} q_{jl} \hat v_l 
  \right].
   \label{eq:RAND.15}
\end{eqnarray}
We dropped the higher-order terms that do not play a role in the limit
$\eps \to 0$, and replaced the equality with the approximate sign. To
simplify our notation, we omit henceforth all the arguments in the
equations, except those of $\mu$.  The arguments will be spelled out
only in definitions.

The coupling matrix in (\ref{eq:RAND.14})-(\ref{eq:RAND.15}) is given by 
\begin{equation}
\label{eq:RAND.17}
  q_{jl} = 2 \int_0^{{\mathcal D}} d \eta \, \phi_j(\eta)
\phi_l''(\eta) = -2 \left(\frac{\pi}{{{\mathcal D}}}\right)^2 \Big(j -
\frac{1}{2}\Big)^2 \delta_{jl}.
\end{equation}
It takes this simple diagonal form because we assumed a homogeneous
background speed $c_o$. If we had a variable speed $c(y)$, the
matrix $\{q_{jl}\}$ would not be diagonal, and the modes with $j \ne
l$ would be coupled. However, the results of the asymptotic analysis
below would still hold, because the coupling would become negligible
in the limit $\eps \to 0$ considered in section \ref{sect:LR}, due to
rapid phases arising in the right hand sides of (\ref{eq:RAND.14}),
(\ref{eq:RAND.15}).

The equations for the evanescent components are obtained similarly,
\begin{equation}
\label{eq:RAND.16}
  \left(\partial_z^2 + \partial_x^2 - \beta_j^2\right) \hat v_j \approx 
  \eps^{3/2} \mu(\eps x,\eps z) \, q_{jj} \hat v_j, 
\end{equation}
and they are augmented with the decay conditions $\hat v_j(\om,x,z) \to
0$ as $\sqrt{x^2+z^2}\to \infty$, for all $j \ge N+1$.

\section{The limit process}
\label{sect:LR}
We characterize next the wave field in the asymptotic limit $\eps \to
0$. We begin with the paraxial long range scaling that gives
significant net scattering, and then take the limit. The scaling has
already been described at the beginning of section \ref{sect:RANDWG}.
\subsection{Asymptotic scaling}
\label{sect:LR.1}
We obtain from (\ref{eq:RAND.14})-(\ref{eq:RAND.17}) that the
propagating mode amplitudes satisfy the block diagonal system of
partial differential equations
\begin{equation}
\left( \begin{array}{cc}
 2 i \beta_j \partial_z + 
 \partial_x^2 & e^{-2 i \beta_j z} \partial_x^2 \\
 e^{2 i \beta_j z} \partial_x^2 & -2 i \beta_j \partial_z + 
    \partial_x^2
\end{array} \right) \left( \begin{array}{c}
a_j \\ b_j \end{array} \right) \approx 
\eps^{3/2} q_{jj} \mu(\eps x, \eps z) \left( \begin{array}{cc}
1 & e^{-2 i \beta_j z}  \\
 e^{2 i \beta_j z} & 1 
\end{array} \right) \left( \begin{array}{c}
a_j \\ b_j \end{array} \right), 
\label{eq:LR.1}
\end{equation}
for $j = 1, \ldots, N$. Again, the approximate sign means equal to
leading order.

Because the right hand side in (\ref{eq:LR.1}) is small, of order
$\eps^{3/2}$, and has zero statistical expectation, it follows from
\cite[Chapter 6]{book07} that there is no net scattering effect until
we reach ranges of order $\eps^{-2}$. Thus, we let
\begin{equation}
z = Z/\eps^2, 
\label{eq:LR.2}
\end{equation}
with scaled range $Z$ independent of $\eps$. The source directivity
in the range direction suggests observing the wavefield on a cross-range 
scale that is smaller than that in range. We choose it as 
\begin{equation}
x = X/\eps,
\label{eq:LR.3}
\end{equation} 
with scaled cross-range $X$ independent of $\eps$, to balance
the two terms in the paraxial operators in
(\ref{eq:LR.1}).  

Our goal is to characterize the $\eps \to 0$ limit of the mode
amplitudes in the paraxial long range scaling regime
(\ref{eq:LR.2})-(\ref{eq:LR.3}). We denote them by
\begin{equation}
  a_j^\eps(\om,X,Z) = a_j\Big(\om,\frac{X}{\eps},\frac{Z}{\eps^2}\Big) 
  \quad \mbox{and} \quad 
  b_j^{\eps} 
  (\om,X,Z) =  
  b_j  \Big(\om,\frac{X}{\eps},\frac{Z}{\eps^2}\Big),
\label{eq:LR.4}
\end{equation}
and obtain from (\ref{eq:LR.1})-(\ref{eq:LR.3}) that they satisfy the
scaled system
\begin{equation}
\left( \begin{array}{cc}
 2 i \beta_j \partial_Z + 
 \partial_X^2 & e^{-2 i \beta_j Z/\eps^2} \partial_X^2 \\
 e^{2 i \beta_j Z/\eps^2} \partial_X^2 & -2 i \beta_j \partial_Z + 
    \partial_X^2
  \end{array} \right) \left( \begin{array}{c}
    a_j^\eps \\ b_j^\eps \end{array} \right) \approx 
\frac{1}{\eps^{1/2}} \mu \Big(X,  \frac{Z}{\eps} \Big) q_{jj} \left( 
  \begin{array}{cc}
1 & e^{-2 i \beta_j Z/\eps^2}  \\
 e^{2 i \beta_j Z/\eps^2} & 1 
\end{array} \right) \left( \begin{array}{c}
a_j^\eps \\ b_j^\eps \end{array} \right), 
\label{eq:LR.5}
\end{equation}
for $j = 1, \ldots, N$, with initial conditions
\begin{equation}
  a_j^\eps(\om,X,0) = a_{j,{\rm ini}}\left(\om,X \right) := \frac{1}{2
  i \beta_j(\omega)} \widehat{F}_j(\omega, X) ,
\label{eq:LR.6}
\end{equation}
and end conditions
\begin{equation}
b_j^\eps(\om,X,L) = 0. 
\label{eq:LR.7}
\end{equation}

\subsection{The random propagator}
\label{sect:LR.1.2}
Let us rewrite (\ref{eq:LR.5}) in terms of
the random propagator matrix ${\bf P}^\eps(\om,X,X',Z) \in
\mathbb{C}^{2 N \times 2 N}$,  the solution of the initial value
problem
\begin{eqnarray}
\partial_Z {\bf P}^\eps(\om,X,X',Z) &=&\left[\frac{1}{\eps^{1/2}} \mu
\Big(X, \frac{Z}{\eps} \Big) {\bf H}\Big(\om,X,\frac{Z}{\eps^2}\Big) +
{\bf G}\Big(\om,X,\frac{Z}{\eps^2}\Big)\right]{\bf
P}^\eps(\om,X,X',Z), \qquad Z > 0, \nonumber \\ {\bf
P}^\eps(\om,X,X',0) &=& \delta(X-X'){\bf I}.
\label{eq:LR.8}
\end{eqnarray}
Here ${\bf I}$ is the $2N \times 2N$ identity matrix, $\delta(X)$
is the Dirac delta distribution in $X$, and ${\bf G}$ and
${\bf H}$ are matrices with entries given by partial differential
operators in $X$, with deterministic coefficients. We can define them
from (\ref{eq:LR.5}) once we note that the solution
\begin{equation}
{\itbf a}^\eps(\om,X,Z) = \left( \begin{array}{c} a_1^\eps (\om,X,Z)\\
\vdots \\ a_N^\eps(\om,X,Z) \end{array}\right), \qquad 
{\itbf b}^\eps(\om,X,Z) =\left( \begin{array}{c} b_1^\eps(\om,X,Z) \\
\vdots \\ b_N^\eps(\om,X,Z) \end{array}\right)
\label{eq:LR.9}
\end{equation}
follows from 
\begin{equation}
\left( \begin{array}{c} 
{\itbf a}^\eps(\om,X,Z) \\
{\itbf b}^\eps(\om,X,Z)
\end{array} \right) = \int d X' \, {\bf P}^\eps(\om,X,X',Z) 
\left( \begin{array}{c}
    {\itbf a}^\eps(\om,X',0) \\{\itbf b}^\eps(\om,X',0) 
  \end{array} \right).
\label{eq:LR.10}
\end{equation}
Here ${\itbf b}^\eps(\om,X',0)$ is the vector of backward going
amplitudes at the beginning of the randomly perturbed section of the
waveguide, and it can be eliminated using the boundary identity
\begin{equation}
\left( \begin{array}{c} 
{\itbf a}^\eps(\om,X,L) \\
{\bf 0}
\end{array} \right) = \int d X' \, {\bf P}^\eps(\om,X,X',L) 
\left( \begin{array}{c} {\itbf a}^\eps(\om,X',0) \\{\itbf
  b}^\eps(\om,X',0) \end{array} \right).
\label{eq:LR.10p}
\end{equation}
The initial conditions ${\itbf a}^\eps(\om,X',0)$ are given in
(\ref{eq:LR.6}).

We obtain from (\ref{eq:LR.5}) that ${\bf H}$ and ${\bf G}$ have the
block form
\begin{equation}
{\bf H} = \left( \begin{array}{cc} {\bf H}^{a} & {\bf H}^{b} \\
\overline{{\bf H}^{b}} & \overline{{\bf H}^{a}} \end{array} \right),
\qquad {\bf G} = \left( \begin{array}{cc} {\bf G}^{a} & {\bf G}^{b} \\
\overline{{\bf G}^{b}} & \overline{{\bf G}^{a}} \end{array} \right),
\label{eq:LR.11}
\end{equation}
where the bar denotes complex conjugation. The blocks are diagonal,
with entries
\begin{equation}
  {\bf H}_{jl}^{a} = -\frac{i \, \delta_{jl}\, q_{jj}}{2 \beta_j},
  \qquad {\bf H}_{jl}^{b} = -\frac{i \, \delta_{jl} \, q_{jj}}{2
  \beta_j} e^{-2 i \beta_j Z/\eps^2} \label{eq:LR.12a}
\end{equation}
and 
\begin{equation}
  {\bf G}_{jl}^{a} = \frac{i\, \delta_{jl}}{2 \beta_j} \partial_X^2,
  \qquad {\bf G}_{jl}^{b} = \frac{i \, \delta_{jl}}{2 \beta_j} e^{-2 i
  \beta_j Z/\eps^2} \partial_X^2, \qquad j, l = 1, \ldots, N.
\label{eq:LR.12b}
\end{equation}
 The entries of the diagonal blocks depend only on the mode indices
and the frequency, via $\beta_j(\om)$.  The entries of the
off-diagonal blocks are rapidly oscillating, due to the large phases
proportional to $Z/\eps^2$.

The symmetry relations satisfied by the blocks in ${\bf H}$ and
${\bf G}$ imply that the propagator has the  form
\begin{equation}
{\bf P}^\eps(\om,X,X',Z) = \left( \begin{array}{cc}
{\bf T}^\eps(\om,X,X',Z) & \overline{{\bf R}^\eps(\om,X,X',Z)} \\
{\bf R}^\eps(\om,X,X',Z) & \overline{{\bf T}^\eps(\om,X,X',Z)}
\end{array} \right),
\label{eq:LR.13}
\end{equation}
with $N \times N$ complex, diagonal blocks ${\bf T}^\eps$ and $
{\bf R}^\eps$.  

\subsection{The diffusion limit}
\label{sect:LR.2}
The limit of ${\bf P}^\eps$ as $\eps \to 0$ is a multi-dimensional
Markov diffusion process, with entries satisfying a system of
It\^{o}-Schr\"{o}dinger equations. This follows from the diffusion
approximation theorem \cite{kohler74,kohler75}, see also \cite[Chapter
6]{book07}, applied to system (\ref{eq:LR.8}).

When computing the generator of the limit process, we obtain that due
to the fast phases in the off-diagonal blocks of ${\bf H}$ and ${\bf
G}$, the forward and backward going amplitudes decouple as $\eps \to
0$. This implies that there is no backscattered field in the limit,
because the backward going amplitudes ${\itbf b}^\eps$ are set to zero at $Z =
L$. Equation (\ref{eq:LR.10}) simplifies as
\begin{equation}
{\itbf a}^\eps(\om,X,Z) = \int d X'\, {\bf T}^\eps(\om,X,X',Z) 
{\itbf a}^\eps(\om,X',0),
\end{equation} 
where the initial conditions ${\itbf a}^\eps(\om,X',0)$ are given in
(\ref{eq:LR.6}).
We call the complex diagonal matrix 
\begin{equation}
{\bf T}^\eps(\om,X,X',Z) = \diag \left( \cT_1^\eps(\om,X,X',Z), \ldots, 
\cT_N^\eps(\om,X,X',Z)\right)
\label{eq:TRMAT}
\end{equation}
the transfer process, because it gives the amplitudes of the forward
going modes at positive ranges $Z$, in terms of the initial conditions
at $Z = 0$.
The limit transfer process is described in the next proposition.
It follows straight from \cite{kohler74,kohler75}.

\vspace{0.05in}
\begin{proposition}
  \label{prop.1} As $\eps \to 0$, ${\bf T}^\eps(\om,X,X',Z)$ converges
  weakly and in distribution to the diffusion Markov process
  ${\bf T}(\om,X,X',Z)$.  This process is complex and diagonal matrix
  valued, with entries $\cT_j(\om,X,X',Z)$ solving the
  It\^{o}-Schr\"{o}dinger equations
\begin{equation}
  d\cT_j(\om,X,X', Z) = \left[\frac{i}{2 \beta_j(\om)} \partial_X^2 -
  \frac{q_{jj}^2 R_o(0) }{8 \beta_j^2(\om)}\right] \cT_j(\om ,X, X',Z)
  dZ + \frac{i \, q_{jj}}{2 \beta_j(\om)} \cT_j(\om, X, X',Z) d
  \cB(X,Z),
\label{eq:LR.14}
\end{equation}
for $Z > 0$, and initial conditions
\begin{equation}
\cT_j(\om,X, X',0) = \delta(X-X'), \qquad j = 1, \ldots, N.
            \label{eq:LR.15}
\end{equation}
Equations (\ref{eq:LR.14}) are uncoupled, but they are driven by the
same Brownian field $\cB(X,Z)$, satisfying
\begin{equation}
\EE\left[ \cB(X,Z) \right] = 0, \quad \EE \left[ \cB(X,Z) \cB(X',Z')
\right] = \min\{Z,Z'\} R_o(X-X'),
\label{eq:LR.16}
\end{equation}
with $R_o$ defined in (\ref{eq:RAND.3p}). Thus, the transfer
coefficients $\cT_j$ are statistically correlated.
\end{proposition}

\vspace{0.05in}
The weak convergence in distribution means that we can calculate the
limit $\eps \to 0$ of statistical moments of ${\bf T}^\eps$,
smoothed by integration over $X$ against the initial conditions, using
the Markov diffusion defined by (\ref{eq:LR.14})-(\ref{eq:LR.15}). In
applications we have a fixed $\eps \ll 1$, and we use Proposition
\ref{prop.1} to approximate the statistical moments of the amplitudes
of the forward going waveguide modes.

When comparing the It\^{o}-Schr\"{o}dinger equations (\ref{eq:LR.14})
to the deterministic Schr\"{o}dinger equations (\ref{eq:HOM.16})
satisfied by the amplitudes in the ideal waveguides, we see that the
random boundary scattering effect amounts to a net diffusion, as
described by the last two terms in (\ref{eq:LR.14}).  We show next how
this leads to loss of coherence of the waves, that is to exponential
decay in range of the mean field. We also study the propagation of
energy of the modes and quantify the decorrelation properties of the
random fluctuations of their amplitudes.

\section{Statistics of the wave field}
\label{sect:STAT}
We begin in section \ref{sect:STAT.1} with the analysis of the
coherent field. Explicitly, we estimate the mean forward going mode
amplitudes in the paraxial long range regime.  Traditional imaging
methods rely on these being large with respect to their random
fluctuations. However, this is not the case, because
$\EE\left[a_j^\eps(\om,X,Z)\right]$ decay exponentially with $Z$, at
rates that increase monotonically with mode indices $j$.  The second
moments of the amplitudes do not decay, but there is decorrelation
over the modes and the frequency and cross-range offsets, as shown in
sections \ref{sect:STAT.2} and \ref{sect:STAT.3}.  Understanding these
decorrelations is key to designing time reversal and imaging methods
that are robust at low SNR.  Robust means that wave focusing in time
reversal or imaging is essentially independent of the realization of
the random boundary fluctuations, it is statistically stable. Low SNR
means that the coherent (mean) field, the ``signal'', is faint with
respect to its random fluctuations, the ``noise''.
\subsection{The coherent field}
\label{sect:STAT.1}
The mean modal amplitudes are 
\begin{eqnarray}
  \EE\left[ a_j^\eps(\om,X, Z)\right] &\approx&  \int d X' \, \EE \left[ 
    \cT_j(\om,X,X',Z) \right]a_{j,{\rm ini}}\left( \om, X' \right),
\label{eq:STAT.1}
\end{eqnarray}
with mean transfer matrix satisfying the partial differential equation
\begin{equation}
  \partial_Z \EE\left[ \cT_j(\om,X,X',Z)\right] = 
\left[\frac{i}{2 \beta_j(\om)} 
    \partial^2_X - \frac{1}{\cS_j(\om)} \right] 
  \EE\left[ \cT_j(\om,X,X',Z)\right], \qquad Z > 0, 
\label{eq:STAT.1p}
\end{equation}
with mode-dependent damping coefficients 
\begin{equation}
  \cS_j(\om) = \frac{8 \beta_j^2(\om)}{q_{jj}^2 R_o(0)} = 
  \frac{2 {{\mathcal D}}^2}{\sigma^2 
    \pi^2 \ell} \left[ \frac{(N+\alpha(\om)-1/2)^2-(j-1/2)^2}{(j-1/2)^4}\right], 
\label{eq:STAT.3}
\end{equation}
with units of length.  Here we used definitions (\ref{eq:HOM.13}),
(\ref{eq:DefAlpha}) and (\ref{eq:RAND.3}), and obtained equation
(\ref{eq:STAT.3}) by taking expectations in (\ref{eq:LR.14}).  Its
solution is given by
\begin{equation}
  \EE \left[ 
    \cT_j(\om,X,X',Z) \right] = \sqrt{\frac{\beta_j(\om)}{2 \pi i Z}} 
  \exp \left[ -\frac{Z}{\cS_j(\om)} + \frac{i \beta_j(\om) (X-X')^2}{2 Z}
\right], 
\label{eq:STAT.2}
\end{equation}
and the mean modal amplitudes are obtained from equations (\ref{eq:STAT.1})
and (\ref{eq:LR.6})
\begin{eqnarray}
  \EE\left[ a_j^\eps(\om,X, Z)\right] &\approx& - \frac{1}{2}
  \sqrt{\frac{i}{2 \pi \beta_j(\om) Z}} \int d X' \, \hat F_j(\om,X')
  \exp \left[ -\frac{Z}{\cS_j(\om)} + \frac{i \beta_j(\om)(X- X')^2}{2
  Z}\right] \nonumber \\ &=& a_{j,o}\left(\om, X, Z \right) \exp
  \left[ -\frac{Z}{\cS_j(\om)} \right],
  \label{eq:STAT.4}
\end{eqnarray}
with $a_{j,o}$ the solution of the paraxial wave equation
(\ref{eq:HOM.16}-\ref{eq:HOM.16b}) in the ideal waveguide.

The mean wave field follows from (\ref{eq:RAND.10}), after neglecting
the evanescent part,
\begin{equation}
  \EE\left[ \hat P\Big(\om,\frac{X}{\eps},\eta,\frac{Z}{\eps^2}\Big)
  \right] \approx \sum_{j=1}^N \phi_j(
 \eta) 
  a_{j,o}\left(\om, X , {Z} \right) 
  \exp \left[ -\frac{Z}{\cS_j(\om)} + i 
    \beta_j(\om) \frac{Z}{\eps^2}  \right].
\label{eq:STAT.5}
\end{equation}
It is different than the field in the ideal waveguides
\begin{equation}
  \hat p_o\Big(\om,\frac{X}{\eps},\eta,\frac{Z}{\eps^2} \Big)
 \approx \sum_{j=1}^N \phi_j( \eta) 
  a_{j,o}\left(\om, X , Z \right) 
  \exp \left[ i 
    \beta_j(\om) \frac{Z}{\eps^2}  \right],
\label{eq:STAT.5id}
\end{equation}
because of the exponential decay of the mean mode amplitudes, on 
range scales $\cS_j(\om)$. 

\subsection{High-frequency and low-SNR regime}
\label{sect:SNR}
We call the length scales $\cS_j(\om)$ the \emph{mode-dependent
  scattering mean free paths}, because they give the range
over which the modes become essentially incoherent, with low SNR,
\begin{equation}
  {\rm SNR}_{j,\om} = \frac{\left|\EE\left[ a_j^\eps(\om,X,Z)\right]\right|}
  {\sqrt{ \EE
      \left[ | a_j^\eps(\om,X,Z)|^2\right] - \left| \EE\left[
          a_j^\eps(\om,X,Z)\right]\right|^2}} \sim 
  \exp\left[ -\frac{Z}{\cS_j(\om)} \right] \ll 1, \quad \mbox{if} ~ ~ 
  Z \gg \cS_j(\om).
\label{eq:STAT.6}
\end{equation}
The second moments $\EE \left[ | a_j^\eps(\om,X,Z)|^2\right]$ are
calculated in the next section, and they do not decay with range. This 
is why equation (\ref{eq:STAT.6}) holds.

The scattering mean free paths decrease monotonically with mode
indices $j$, as shown in (\ref{eq:STAT.3}). The first mode encounters
less often the random boundary, and has the longest scattering
mean free path
\begin{equation}
  \cS_1(\om) = \frac{32 {{\mathcal D}}^2}{\sigma^2 \pi^2 \ell} \left[
    (N+\alpha(\om)-1/2)^2 - 1/4\right] \approx \frac{32 {{\mathcal D}}^2 N^2}
{\sigma^2 \pi^2 \ell}.
\label{eq:S1}
\end{equation}
The highest indexed mode scatters most frequently at the boundary,
and its scattering mean free path 
\begin{equation}
  \cS_N(\om) = \frac{2 {{\mathcal D}}^2}{\sigma^2 \pi^2 \ell}
      \frac{\alpha(\om)\left( 2N + \alpha(\om)-1\right)}{(N-1/2)^4}
      \approx \frac{\alpha(\om)}{8} \frac{\cS_1(\om)}{N^5}  
\label{eq:SN}
\end{equation}
is much smaller than $\cS_1(\om)$, when $N$ is large.
To be complete, we also have
$$
  \cS_j(\om) \approx   \cS_1(\om) \frac{1-s^4}{s^4} \frac{1}{N^4} ,
   \quad \quad \mbox{ if } j =\lfloor sN \rfloor, \quad s\in(0,1),
$$
and
$$
  \cS_j(\om) \approx  \cS_1(\om) \frac{1}{(2j-1)^4} ,
   \quad \quad \mbox{ if } j =o(N).
$$

Our analysis of time reversal and imaging is carried in a
high-frequency regime, with waveguide depth $\mathcal{D}$ much larger
than the central wavelength $\la_o$ or, equivalently, with $N \gg 1$.
We also assume a low-SNR regime, with scaled range $Z$ exceeding the
scattering mean free path of all the modes, so that none of the
amplitudes $a_j$ are coherent.  This is the most challenging case for
sensor array imaging, because the wave field measured at the sensors
is essentially just noise.  We model the low-SNR regime using the
dimensionless large parameter
\begin{equation}
\gamma = \frac{Z}{\cS_1(\om_0)} \gg 1,
\label{eq:gamma}
\end{equation}
and observe from (\ref{eq:STAT.3}) that 
\begin{equation}
  \frac{Z}{\cS_j(\om_0)} 
      \ge \gamma \gg 1, 
\quad \mbox{for all } j = 1, \ldots, N.
\label{eq:gamma1}
\end{equation}

\subsection{The second moments}
\label{sect:STAT.2}
The quantification of SNR and the analysis of time reversal and
imaging involves the second moments of the mode amplitudes.  Recall
that
\begin{equation}
  a_j^\eps(\om,X,Z) \approx \int dX'  \cT_j^\eps 
  \left(\om,X,X',Z\right)
  a_{j,{\rm ini}}(\om,X' ) 
\end{equation}
with $\cT_j^\eps$ the entries of the diagonal transfer matrix ${\bf
T}^\eps$. To calculate the second moments, we need to estimate
$\EE\left[ \cT_j^\eps \overline{\cT_l^\eps}\right]$.  The equations
for $\cT_j^\eps(\om_1,X_1,X_1',Z)
\overline{\cT_l^\eps(\om_2,X_2,X_2',Z)}$ follow from the forward
scattering approximation of (\ref{eq:LR.8}),
\begin{eqnarray}
  \partial_Z \cT_j^\eps \overline{\cT_l^\eps} &\approx&
  \left[\frac{i}{2 \beta_j(\om_1)}\partial_{X_1}^2 - \frac{i}{2
      \beta_l(\om_2)}\partial_{X_2}^2\right] \cT_j^\eps
  \overline{\cT_l^\eps} \nonumber \\ && - \frac{i}{2 \eps^{1/2}}\left[
    \frac{q_{jj}\,\mu(X_1,Z/\eps)}{ \beta_j(\om_1)} -
    \frac{q_{ll}\,\mu(X_2,Z/\eps)}{ \beta_l(\om_2)}\right]\cT_j^\eps
  \overline{\cT_l^\eps},
\end{eqnarray}
for $Z>0$, with initial condition
\begin{equation}
  \cT_j^\eps(\om_1,X_1,X_1',0) \overline{\cT_l^\eps(\om_2,X_2,X_2',0)} = 
\delta(X_1-X_1')\delta(X_2-X_2').
\end{equation}
Their statistical distribution is characterized in the limit $\eps \to
0$ by the diffusion approximation theorem \cite{kohler74,kohler75},
see also \cite[Chapter 6]{book07}.  It is the distribution of
$\cT_j(\om_1,X_1,Z) \overline{\cT_l}(\om_2,X_2,Z)$, with $\cT_j$ the
limit transfer coefficients in Proposition \ref{prop.1}. This gives
the approximate relation
\begin{eqnarray}
  \EE\left[ a_j^\eps(\om_1,X_1,Z) \overline{a_l^\eps(\om_2, X_2,Z)}\right] 
 & \approx& \int dX_1' \int dX_2' \, 
  a_{j,{\rm ini}}\left(\om_1, X_1' \right)   \overline{ 
    a_{l,{\rm ini}}\left(\om_2,X_2'\right)}\nonumber \\ 
 &&\times  \EE \left[ \cT_j(\om_1,X_1,X_1',Z) \overline{
      \cT_l(\om_2, X_2,X_2',Z)}\right].\label{eq:STAT.7}
\end{eqnarray}
The calculation of $\EE\left[\cT_j \overline{\cT_l}\right]$ is given
in appendix \ref{ap:SM}. We summarize the results in
Propositions \ref{prop.2}-\ref{prop.3}.

\subsubsection{The single mode and frequency moments}
It is easier to calculate the diagonal moments, with $j = l$, and the
same frequency $\om_1 = \om_2 = \om$. We have the following result
proved in appendix \ref{ap:SM}.

\vspace{0.05in}
\begin{proposition}
  \label{prop.2} For all $j = 1, \ldots, N$, and all the frequencies $\om 
\in [\om_0-\pi B,\om_0+\pi B]$,
\begin{eqnarray}
  \EE\left[ \cT_j(\om,X_1,X_1',Z)
  \overline{\cT_j(\om,X_2,X_2',Z)}\right] = \frac{\beta_j(\om)}{2 \pi
  Z} \exp \left\{ \frac{i \beta_j(\om) [(X_1-X_1')^2-(X_2-X_2')^2]}{2
  Z} \right. \nonumber \\ \left.  -\frac{2 Z}{\cS_j(\om)}\int_0^1 d s \,
  C_o\big[ (X_1-X_2) s + (X_1'-X_2')(1-s)\big] \right\},
\label{eq:STAT.13}
\end{eqnarray}
with kernel $C_o$  defined by 
\begin{equation}
C_o(X) = 1- \frac{R_o(X)}{R_o(0)}.
\label{eq:STAT.11}
\end{equation}
\end{proposition}

\vspace{0.05in}
The general second moment formula does not have an explicit form in
arbitrary regimes. But it can be approximated in the 
low-SNR regime (\ref{eq:gamma}). The expression
(\ref{eq:STAT.13}) also simplifies in that regime, 
as stated in the following proposition, which we prove below.

\vspace{0.05in}
\begin{proposition}
\label{prop.2ap}
In the low-SNR regime (\ref{eq:gamma}), and under the assumption $X_1'
= X_2' = X'$, the right hand side in (\ref{eq:STAT.13}) is essentially
zero, unless
\begin{equation}
  \frac{|X_1-X_2|}{\ell} \lesssim \sqrt{\frac{3 \, \cS_j(\om)}
{\gamma \, \cS_1(\om)}} \ll 1, 
\label{eq:REM.1}
\end{equation}
and the moment formula simplifies to 
\begin{equation} \EE\left[
    \cT_j(\om,X_1,X',Z) \overline{\cT_j(\om,X_2,X',Z)}\right] \approx
  \frac{\beta_j}{2 \pi Z} \exp \left[ \frac{i
  \beta_j[(X_1-X')^2-(X_2-X')^2]}{2 Z}- \frac{(X_1-X_2)^2}{2
  X_{d,j}^2(\om)} \right],
\label{eq:STAT.14}
\end{equation}
with 
\begin{equation}
  X_{d,j}(\om) = \ell \sqrt{ \frac{3 \, \cS_j(\om)}{2 Z}} =
 \ell \sqrt{ \frac{3 \, \cS_j(\om)}{2 \gamma \, \cS_1(\om)}} \ll \ell.
\label{eq:STAT.15}
\end{equation}
If the initial points $X_1'$ and $X_2'$ are different, but still close
enough to satisfy
\begin{equation}
  \frac{|X_1'-X_2'|}{\ell}  \ll 1, 
\label{eq:REM.1p}
\end{equation}
the moment formula becomes 
\begin{eqnarray} 
\EE\left[ \cT_j(\om,X_1,X_1',Z)
    \overline{\cT_j(\om,X_2,X_2',Z)}\right] \approx \frac{\beta_j}{2
    \pi Z} \exp \left[ \frac{i \beta_j[(X_1-X_2')^2-(X_2-X_2')^2]}{2
    Z}\right] \nonumber \\ \times \, \exp \left[ -\frac{(X_1-X_2)^2 +
    (X_1'-X_2')^2+(X_1-X_2)(X_1'-X_2')}{2 X_{d,j}^2(\om)} \right].
\label{eq:STAT.14p}
\end{eqnarray}
\end{proposition}
\vspace{0.05in}
\begin{proof}
We see from definitions (\ref{eq:RAND.3p}) and (\ref{eq:STAT.11}) that
$C_o(X) \approx 1$ for $|X| \gg \ell$. Therefore,
\[
\int_0^1 d s \, C_o\big[ (X_1-X_2) s \big] \approx 1 \quad \mbox{if} ~ ~
|X_1 - X_2| \gg \ell,
\]
and the right hand side in (\ref{eq:STAT.13}) becomes negligible, of
order $ O\left(e^{-2 Z/\cS_j}\right) \ll 1.  $
In the case $|X_1 - X_2| \sim \ell$ we obtain similarly that the
damping term is of order $Z$, and the right hand
side in (\ref{eq:STAT.13}) is exponentially small.
It is only when $|X_1-X_2| \ll \ell$ that the moment does not vanish.
Then, we can approximate the kernel $C_o$ in the integral with its first
nonzero term in the Taylor expansion around zero, using the relations
\begin{equation}
  C_o(0) = 0, \quad C_o'(0) = 0, ~ ~ \mbox{and} ~ 
C_o''(0) = -\frac{R_o''(0)}{R_o(0)} = \frac{1}{\ell^2},
\label{eq:STAT.12}
\end{equation}
that follow from (\ref{eq:RAND.3})-(\ref{eq:RAND.4}). We have 
\[
\int_0^1 d s \, C_o \big[  (X_1-X_2) s \big] \approx 
\frac{|X_1-X_2|^2 }{6 \ell^2}, \quad 
\mbox{if} ~ ~ |X_1 - X_2| \gg \ell,
\]
and the right hand side in (\ref{eq:STAT.13}) is of the order
$\exp\left[-\frac{|X_1-X_2|^2 Z}{3 \ell^2 \cS_j} \right]$. This gives
the condition (\ref{eq:REM.1}), and the simpler moment formula
(\ref{eq:STAT.14}) follows. 

Essentially the same proof applies in the case $X_1'\ne X_2'$, because
we can still expand the integrand in (\ref{eq:STAT.13}) by assumption
(\ref{eq:REM.1p}).
\end{proof}

\subsubsection{The two mode and frequency moments}
The general second moment formula is derived in appendix \ref{ap:SM},
in the low-SNR regime (\ref{eq:gamma}). It has a complicated expression 
that we do not repeat here, but it simplifies for nearby frequencies, as
stated below. 

\vspace{0.05in}
\begin{proposition}
\label{prop.3}
The modes decorrelate under the low-SNR assumption (\ref{eq:gamma})
\begin{equation}
  \EE\left[ \cT_j\left(\om_1,X_1,X_1',Z\right) \overline{\cT_l
      \left(\om_2,X_2,X_2',Z\right)}\right] \approx 0 \quad 
\mbox{if} ~ j \ne l,
\label{eq:2M.1}
\end{equation}
for any two frequencies $\om_1,\om_2$ and cross-ranges $X_1, X_2$.
The modes also decorrelate for frequency offsets that exceed 
\begin{equation}
  \Omega_{d,j}(\om) = \frac{\cS_j(\om)
    \beta_j^2(\om) \ell^2}{Z^2 |\beta'_j(\om)|} = 
  \frac{\beta_j(\om)}{|\beta'_j(\om)|}\frac{\cS_j(\om) 
    \beta_j(\om)\ell^2}{\gamma^2 \cS_1^2(\om)},
\label{eq:2M.3}
\end{equation}
where $\beta'_j(\om)$ is the derivative of $\beta_j(\om)$ with
respect to $\om$.  For much smaller frequency offsets satisfying
\begin{equation}
  |\om_1-\om_2| \ll \Omega_{d,j}\left(\om\right), \quad 
  \om = \frac{\om_1+\om_2}{2},
\end{equation}
the moment formula is 
\begin{eqnarray}
  \EE\left[ \cT_j\left(\om_1,X_1,X_1',Z\right) \overline{\cT_j
      \left(\om_2,X_2,X_2',Z\right)}\right] 
  \approx \frac{\beta_j(\om)}{2 \pi Z} 
  \exp \left\{ \frac{i \left[\beta_j(\om_1) (X_1-X_1')^2 - \beta_j(\om_2)
        (X_2-X_2')^2\right]}{2 Z}  \right. 
  \nonumber \\
  \left. -
    \frac{(X_1-X_2)^2 + (X_1'-X_2')^2 + (X_1-X_2)(X_1'-X_2')}{2 X_{d,j}^2(\om)}
  \right\}.
  \qquad \quad 
\label{eq:2M.4}
\end{eqnarray}
\end{proposition}
\subsection{Decorrelation properties}
\label{sect:STAT.3}
We already stated the decorrelation of the modes in Proposition
\ref{prop.3}. But even for a single mode, we have decorrelation over
cross-range and frequency offsets.

The {\em decoherence length} of mode $j$ is denoted
by $X_{d,j}(\om)$, and it is defined in (\ref{eq:STAT.15}). It is the
length scale over which the second moment at frequency $\om$ decays
with cross-range. It follows from (\ref{eq:STAT.15}) that $X_{d,j}$ is
much smaller than the correlation length, for all the modes, and that
it decreases monotonically with $j$. The first mode has the largest
decorrelation length
\begin{equation}
X_{d,1}(\om) = \ell \sqrt{\frac{3}{2 \gamma}}.
\label{eq:2M.5}
\end{equation}
because it scatters less often at the boundary. 
The decoherence length of the highest mode is much smaller in 
high-frequency regimes with $N \gg 1$,
\begin{equation}
  X_{d,N}(\om) = X_{d,1}(\om) \sqrt{\frac{\cS_N(\om)}{\cS_1(\om)}} 
\approx \frac{\ell}{8} \sqrt{\frac{3 \alpha(\om)}{\gamma}} N^{-5/2}.
\label{eq:2M.6}
\end{equation}

The decorrelation frequency is derived in appendix \ref{ap:2fM}. It is
given by (\ref{eq:2M.3}) or, more explicitly, by
\begin{equation}
  \Omega_{d,j}(\om) \approx \frac{\om \sigma^2 \pi^3}{64 \gamma^2} 
\left( \frac{\ell}{\lambda}\right)^3
  \frac{\left[ \left(N+\alpha(\om)-\frac{1}{2}\right)^2 - 
      \left(j-\frac{1}{2}\right)^2\right]^{5/2}}{N^9 (j-1/2)^4},
\label{eq:2M.10}
\end{equation} 
it is much smaller than $\om$ for all the modes, and it decreases
monotonically with $j$, starting from 
\begin{equation}
  \Omega_{d,1}(\om) \approx \frac{\om \sigma^2 \pi^3(\ell/\lambda)^3}{4 
    \gamma^2N^4} .
\end{equation}

\section{The forward model}
\label{sect:FM}
Let us gather the results and summarize them in the following model of
the pressure field
\begin{equation}
  \hat P\Big(\om,\frac{X}{\eps},\eta,\frac{Z}{\eps^2}\Big) \sim
  \sum_{j=1}^N \frac{\phi_j(\eta)}{2 i \beta_j(\om)} e^{i
  \beta_j(\om)\frac{Z}{\eps^2}} \int d X' \, \cT_j(\om,X,X',Z)
  \int_0^{\mathcal D} d \eta' \phi_j(\eta') \, \hat f (\om,X',\eta') ,
\label{eq:K.1}
\end{equation}
where the symbol $\sim$ stands for approximate, in distribution.  That
is to say, the statistical moments of the random pressure field $\hat
P$ are approximately equal to those of the right handside.  The first
and second moments follow from Propositions \ref{prop.1}-\ref{prop.3}.
In our analysis of time reversal and imaging we take small frequency
offsets, satisfying $|\tom| \ll \Omega_{d,j}(\om)$, so that we can use
the simpler moment formula (\ref{eq:2M.4}). 

The computation of the fourth moments of the transfer coefficients is
quite involved. We estimate in appendix \ref{ap:FM} some of them, for
a particular combination of the mode indices and arguments. These
moments are used in the next sections to show the statistical
stability of the time reversal and coherent interferometric imaging
functions.

We analyze next time reversal and imaging in the low SNR regime, and
assume for convenience that the source (\ref{scaledsource}) has the
separable form
\begin{equation}
  f(t,X,\eta) = \frac{\varphi(t) }{\theta_X \theta_\eta} \rho\left(
    \frac{X-X^\star}{\theta_X},\frac{\eta-\eta^\star}{\theta_\eta}\right),
\label{eq:Ksource1}
\end{equation}
meaning that the same pulse $\varphi(t)$ is emitted from all the
points in the support of the non-negative source density $\rho$. We
scale this support with the dimensionless parameters $\theta_X$ and
$\theta_\eta$, and normalize the source by
\begin{equation}
  \int \frac{d X'}{\theta_X} \int \frac{d \eta'}{\theta_\eta}
  \rho\left(
    \frac{X-X^\star}{\theta_X},\frac{\eta-\eta^\star}{\theta_\eta}\right) = 1.
\label{eq:rrho}
\end{equation}

The coefficients
\begin{equation}
  \hat F_j(\om,X) = \frac{\hat \varphi(\om)}{\theta_X \theta_\eta}
 \int_0^{\mathcal D} d \eta \,
  \phi_j(\eta) \rho\left(
    \frac{X-X^\star}{\theta_X},\frac{\eta-\eta^\star}{\theta_\eta}\right),
  \label{eq:Ksource2}
\end{equation}
are proportional to the Fourier coefficients $\hat \varphi(\om)$ of
the pulse, and are thus supported in the frequency interval $[\om_o -
 B/2, \om_o + B/2].  $ The bandwidth $B$ is small enough so that we can
freeze the number of propagating modes to that at the central
frequency, as explained in section \ref{sect:form.Prop}.  The width of
the pulse $\varphi(t)$ is inverse proportional to $B$, and we
distinguish two regimes: The \emph{broadband} regime with $B \gg
\eps^2 \om_o $, and the \emph{narrowband} regime with $B \le \eps^2
\om_o$. The comparison with $\eps^2$ is because the source is at
range $Z_\cA/\eps^2$ from the array, and the modes arrive at time
intervals of order $1/\eps^2$.  Broadband pulses have smaller support
than these travel times, meaning that we can observe the different
arrivals of the modes, at least in the ideal waveguides.

To analyze the resolution of time reversal and imaging, we study in
detail the case of a source density localized around the point
$(X^\star,\eta^\star,0)$. We say that we study the point spread time
reversal and imaging functions, because the source has small
support.  Note however that it is not a point source. Its support
is quantified by the positive parameters $\theta_X$ and $\theta_\eta$
which are small, but independent of $\eps$. 

\section{Time reversal}
\label{sect:TR}

Let us denote by $D(t,X,\eta)$ the pressure field measured in a time
window $\psi(t/T^\eps)$ at an array $\cA$, with aperture modeled by the
indicator function
\begin{equation}
1_{\cA}(X,\eta) = 1_{\cA_X}(X) 1_{\cA_\eta}(\eta),
\label{eq:TR.1}
\end{equation}
at range $ z_\cA = {Z_\cA}/{\eps^2}.  $ Here $X$ is the scaled
cross-range in the array, related to the cross-range $x$ by $x =
X/\eps$, and $\cA_X\subset \mathbb{R}$ and $\cA_\eta\subset
[0,{{\mathcal D}}]$ are intervals in $X$ and $\eta$.  The window
$\psi$ is a function of dimensionless arguments, of support of order
one, and $T^\eps$ denotes the length of time of the measurements.
Because the waves travel distances of order $\eps^{-2}$, we scale
$T^\eps$ as $ T^\eps = T/\eps^2 $, with $T$ of order one.

In time reversal, the array takes the recorded field $D(t,X,\eta)$,
time reverses it and emits $D(T^\eps-t,X,\eta)$ back in the medium. We
study in this section the resolution of the refocusing of the waves at
the source, in the high-frequency and low-SNR regime described in
section \ref{sect:SNR}. Because we have a random waveguide, the
resolution analysis includes that of statistical stability, given in 
section \ref{sect:TRStatStab}.

\subsection{Mathematical model of  time reversal}
\label{sect:TR.1}
We have in our notation
\begin{equation}
  D(t,X,\eta) = 1_{\cA}(X,\eta) \psi\Big(\frac{t}{T^\eps}\Big)
  P\Big(t,x= \frac{X}{\eps},\eta, z_\cA=\frac{Z_\cA}{\eps^2} \Big),
\label{eq:TR.4}
\end{equation}
with mathematical model following from (\ref{eq:K.1}),
\begin{eqnarray}
  D(t,X,\eta) \approx 1_{\cA}(X,\eta) \psi\Big(\frac{t}{T^\eps}\Big)
  \sum_{j=1}^N \phi_j(\eta) \int d \om \, \frac{e^{i \beta_j(\om)
      \frac{Z_\cA}{\eps^2} -i \om t}}{2 i \beta_j(\om)}
  \int  d X' \, \hat F_j(\om,X')
  \cT_j(\om,X,X',Z_\cA).
\label{eq:TR.5}
\end{eqnarray}
The time reversed field
\begin{equation}
\DTR(t,X,\eta) = D(T^\eps-t,X,\eta)
\label{eq:TR.6}
\end{equation}
has Fourier transform
\begin{equation}
  \hDTR(\om,X,\eta) = e^{i \om T^\eps} \overline{\hat D(\om,X,\eta)}
\label{eq:TR.7}
\end{equation}
with
\begin{eqnarray}
  \hat D(\om,X,\eta) \approx 1_{\cA}(X,\eta)
  \sum_{j=1}^N \phi_j(\eta)\int \frac{d u}{2
    \pi} \, \hat \psi(u)
  \frac{e^{i \beta_j(\om-{\eps^2 u}/{T})\frac{Z_\cA}{\eps^2} }
}{2 i \beta_j(\om-\eps^2 u/T)}
  \nonumber \\
  \times   \int  d X' \, \hat F_j \Big(\om-\frac{\eps^2 u}{T},X'\Big)
  \cT_j \Big(\om-\frac{\eps^2 u}{T},X,X',Z_\cA\Big).
\label{eq:TR.8}
\end{eqnarray}
The small frequency shifts $\eps^2 u/T$ are due to the time scaling,
and we can neglect them in the source terms $\hat F_j$ and in the
amplitude factor $1/\beta_j$.

The model of the observed wave field at search locations $
(x^s,\eta^s,z^s) = \left(\frac{X^s}{\eps},\eta^s,\frac{Z^s}{\eps^2}
\right) $ is given by
\begin{eqnarray}
  \cO(t,X^s,\eta^s,Z^s) = \sum_{j=1}^N \phi_j(\eta^s)
  \int \frac{d \om}{2 \pi}
  \,  \frac{\exp \left[ i \beta_j(\om)
      \frac{Z_\cA-Z^s}{\eps^2} -i \om t
    \right]}{2 i \beta_j(\om)} \nonumber \\
  \times\int d X \int d \eta  \, \phi_j(\eta) \, \hDTR(\om,X,\eta)
  \, \cT_j(\om,X,X^s,  Z_\cA-Z^s),
\label{eq:TR.9}
\end{eqnarray}
using reciprocity. Note the similarity with equation (\ref{eq:K.1}),
except that the source is now at the array, which we approximate in
(\ref{eq:TR.9}) as a continuum, instead of a discrete collection of
sensors. This approximation is convenient for the analysis, because
sums over the sensors are replaced by integrals over the $X$ and
$\eta$ apertures, of lengths $|\cA_X|$ and $|\cA_\eta|$.

Using (\ref{eq:TR.7}) in (\ref{eq:TR.9}) and letting
\begin{equation}
  \Gamma_{jl} = \int_0^{{\mathcal D}} d \eta \, 1_{\cA_\eta}(\eta)\,
\phi_j(\eta) \phi_l(\eta)
\label{eq:TR.10},
\end{equation}
we obtain
\begin{eqnarray}
  \cO(t,X^s,\eta^s,Z^s) \approx \int \frac{d X'}{\theta_X} \int
  \frac{d \eta'}{\theta_\eta} \, \rho
  \left(\frac{X'-X^\star}{\theta_X},
    \frac{\eta'-\eta^\star}{\theta_\eta}\right)
  \sum_{j,l=1}^N \Gamma_{jl} \int \frac{d \om}{2 \pi} e^{i \om
    (T^\eps-t)} \overline{\hat \varphi(\om)} ~ \frac{\phi_j(\eta^s)
    \phi_l(\eta')}{4 \beta_j(\om) \beta_l(\om)} ~~ \nonumber \\ \times
  \int \frac{du}{2 \pi} \, \overline{\hat \psi(u)} \int d X \,
  1_{\cA_X}(X) \, \cT_j \Big(\om,X,X^s,Z_\cA-Z^s \Big)
  \overline{\cT_l\Big(\om-\frac{\eps^2 u}{T},X,X',Z_\cA \Big)}
  \nonumber~~~ \\ \times \exp \left[i \beta_j(\om)
    \frac{Z_\cA-Z^s}{\eps^2} - i \beta_l\Big( \om-\frac{\eps^2
      u}{T}\Big) \frac{Z_\cA}{\eps^2} \right]. ~~~~
  \label{eq:TR.12}
\end{eqnarray}
We define the time reversal function by

\begin{eqnarray}
  \ITR(X^s,\eta^s) = \cO(t=T^\eps,X^s,\eta^s,Z^s = 0).
\label{eq:TR.13}
\end{eqnarray}
It models the wave field observed at the time instant $t = T^\eps$, at
the source range $Z_\cA$. This is when and where the refocusing
occurs.

In the case of a source density that is tightly supported around
$(X^\star,\eta^\star)$, we may approximate $\ITR$ by
\begin{equation}
\ITR(X^s,\eta^s) \approx
  \int \frac{d \om}{2 \pi}
\,  \overline{\hat \varphi(\om)}\, \FTR(\om,X^s,\eta^s),
\label{eq:TR.13.1}
\end{equation}
with frequency-dependent kernel (point spread function)
\begin{eqnarray}
  \FTR(\om,X^s,\eta^s) \approx \sum_{j,l=1}^N \Gamma_{jl}
  \frac{\phi_j(\eta^s) \phi_l(\eta^\star)}{4 \beta_j(\om)
  \beta_l(\om)} \int \frac{du}{2 \pi} \, \overline{\hat \psi(u)} \int
  d X \, 1_{\cA_X}(X) \nonumber \\ \times \, \cT_j
  \Big(\om,X,X^s,Z_\cA\Big) \overline{\cT_l\Big(\om-\frac{\eps^2
  u}{T},X,X^\star,Z_\cA \Big)} \nonumber \\ \times \exp \left[i
  \beta_j(\om) \frac{Z_\cA}{\eps^2} - i \beta_l\Big(
  \om-\frac{\eps^2 u}{T}\Big) \frac{Z_\cA}{\eps^2} \right].
  \label{eq:TR.12p}
\end{eqnarray}
Here we used the source normalization (\ref{eq:rrho}).

\subsection{Resolution analysis}
\label{sect:TR.2}
If the time reversal process is statistically stable, then we can
estimate its refocusing resolution by studying the mean of
(\ref{eq:TR.13.1}). We refer to the next section for the analysis of the
statistical stability of $\ITR$.

The mean time reversal function follows from (\ref{eq:2M.4}) and
(\ref{eq:TR.13.1})-(\ref{eq:TR.12p})
\begin{equation}
  \EE\left[\ITR(X^s,\eta^s)\right] =
\int \frac{d \om}{2 \pi} \,
  \overline{\hat \varphi(\om)} \,
  \EE\left[\FTR(\om,X^s,\eta^s)\right],
\label{eq:TR.14}
\end{equation}
with
\begin{eqnarray}
  \EE\left[\FTR(\om,X^s,\eta^s)\right]\approx \frac{|\cA_X|}{8 \pi
  Z_\cA}\sum_{j=1}^N \Gamma_{jj} \frac{\phi_j(\eta^s)
  \phi_j(\eta^\star)}{\beta_j(\om)} \psi \left(\frac{ 
  \beta'_j(\om) Z_\cA}{T}\right)\exp \left[ - \frac{(X^s-X^\star)^2}{2
  X^2_{d,j}(\om)}\right] \nonumber \\ \times \int dX \,
  \frac{1_{\cA_X}(X)}{|\cA_X|} \exp \left[ -\frac{i
  \beta_j(\om)}{Z_\cA}\Big(X-\frac{X^s + X^\star}{2} \Big)
  (X^s-X^\star)\right].
\label{eq:TR.15}
\end{eqnarray}
Moreover, letting
\begin{equation}
\label{eq:DefAX}
\cA_X =\left[-\frac{|A_X|}{2},\frac{|A_X|}{2}\right],
\end{equation}
we obtain after integrating in $X$ that
\begin{eqnarray}
  \EE\left[\FTR(\om,X^s,\eta^s)\right] \approx \frac{|\cA_X|}{8 \pi
Z_\cA}\sum_{j}^N \Gamma_{jj} \frac{\phi_j(\eta^s)
\phi_j(\eta^\star)}{\beta_j(\om)} \psi \left(\frac{ 
\beta'_j(\om) Z_\cA}{T}\right) \exp \left[ - \frac{(X^s-X^\star)^2}{2
X^2_{d,j}(\om)}\right] \nonumber \\ \times \mbox{\rm sinc} \left[
\frac{\beta_j(\om)|\cA_X|}{2 Z_\cA} (X^s-X^\star) \right] \exp \left\{
\frac{i \beta_j(\om)}{2Z_\cA}\left[(X^s)^2 - (X^\star)^2
\right]\right\}.
\label{eq:TR.16}
\end{eqnarray}
Note that
\begin{equation}
\tau_j = \beta'_j(\om) Z_\cA
\end{equation}
are the scaled travel times of the modes, so only those modes that
arrive within the support of the window $\psi$ contribute in
(\ref{eq:TR.16}).

\subsubsection{Cross-range  resolution}
We observe in (\ref{eq:TR.16}) that modes contribute differently to
the focusing in cross-range $X$, with resolution
\begin{equation}
  |X^s-X^\star| \le \Delta_{X,j}(\om):= \min\left\{ X_{d,j}(\om),
\frac{2 \pi Z_\cA}{\beta_j(\om)|\cA_X|} \right\}.
\label{eq:TR.17}
\end{equation}
Recall from (\ref{eq:STAT.15}) and (\ref{eq:2M.5}) that $X_{d,j}$
decreases monotonically with $j$
\begin{equation}
  X_{d,j}(\om) \approx \frac{X_{d,1}(\om)}{4(j-1/2)^2} \frac{\left[
(N+\alpha(\om)-1/2)^2 - (j-1/2)^2\right]^{1/2}}{N}, \qquad
X_{d,1}(\om) = \ell \sqrt{\frac{3}{2 \gamma}},
\label{eq:TR.18}
\end{equation}
whereas
\begin{equation}
\frac{2 \pi Z_\cA}{\beta_j(\om)|\cA_X|} \approx
\frac{2 Z_\cA {{\mathcal D}}}{|\cA_X|} \left[
(N+\alpha(\om)-1/2)^2 - (j-1/2)^2\right]^{-1/2},
\label{eq:TR.19}
\end{equation}
increases with $j$. Thus, in the high-frequency regime with $N \gg 1$,
the cross-range resolution for the high-order modes is determined by
the decorrelation length, even for large apertures. The cross-range
resolution of the first modes may be determined by the aperture, but
only if it is large enough,
\begin{equation}
  |\cA_X| \gtrsim \frac{2 Z_\cA {{\mathcal D}}}{\ell}
  \sqrt{\frac{2 \gamma}{3}} N.
\label{eq:TR.20}
\end{equation}
It may appear at this point that the time reversal process can give
good results even for small apertures $|\cA_X|$. However, we will see
in section \ref{sect:TRStatStab} that large apertures are needed for
statistical stability.

The modes with higher indices give the best cross-range resolution,
but they travel at smaller speed. Thus, the focusing improves when we
increase the recording time, because the array can capture the
late arrivals of the high-order modes (see Figure \ref{fig_parax1}).

\subsubsection{Depth resolution} To study the focusing in $\eta$, we 
evaluate the point spread function at cross-range $X^s = X^\star$. We have
\begin{equation}
  \FTR(\om,X^\star,\eta^s,X^\star,\eta^\star) \approx \sum_{j=1}^{N_T}
\Gamma_{jj} \frac{\phi_j(\eta^s) \phi_j(\eta^\star)}{\beta_j(\om)},
\label{eq:TR.21}
\end{equation}
where $N_T$ is the number of modes with arrival times in the recording
window,
\begin{equation}
\tau_j < T, \qquad \mbox{for} ~ j = 1, 2, \ldots, N_T \le N.
\label{eq:TR.22}
\end{equation}
The coefficients $\Gamma_{jj}$ are given by
\begin{equation}
  \Gamma_{jj} = \int_0^{\mathcal D} d \eta \, 1_{\cA_\eta}(\eta)
  \phi_j^2(\eta) =
  \frac{|\cA_\eta|}{{{\mathcal D}}} + \frac{\eta_{2}}{{\mathcal D}}
  \mbox{sinc}\left[ 2 \pi (j-1/2) \frac{\eta_2}{{{\mathcal D}}}
  \right] - \frac{\eta_{1}}{{\mathcal D}} \mbox{sinc}\left[ 2 \pi (j-1/2)
    \frac{\eta_1}{{{\mathcal D}}}
  \right]  \ge 0,
\label{eq:TR.23}
\end{equation}
for an array in the set $ \cA_\eta = [\eta_1,\eta_2] \subset
[0,{{\mathcal D}}].  $ They satisfy $\Gamma_{jj} = 1$ in the full
aperture case $\cA_\eta = [0,{{\mathcal D}}].$

The sum in (\ref{eq:TR.21}) is maximum at $ \eta^s = \eta^\star$,
because all the terms are positive.  The point spread function is
smaller at other depths, because of cancellations in the sum of the
oscillatory terms. We can make this more explicit
in the high-frequency regime, with $N \gg 1$, if we write
\begin{equation}
\label{eq:HF.1}
{{\mathcal D}} \approx \frac{\pi N}{k}, \qquad
\phi_j(\eta) \approx \sqrt{\frac{2}{{{\mathcal D}}}} \cos \left[ \Big(
    j-\frac{1}{2}\Big)\frac{k\eta}{N}\right],
\end{equation}
and interpret (\ref{eq:TR.21}) as a Riemann sum, which we then
approximate with an integral.

\begin{figure}
   \vspace*{-1.2in}
      \begin{center}
   \includegraphics[width=8cm]{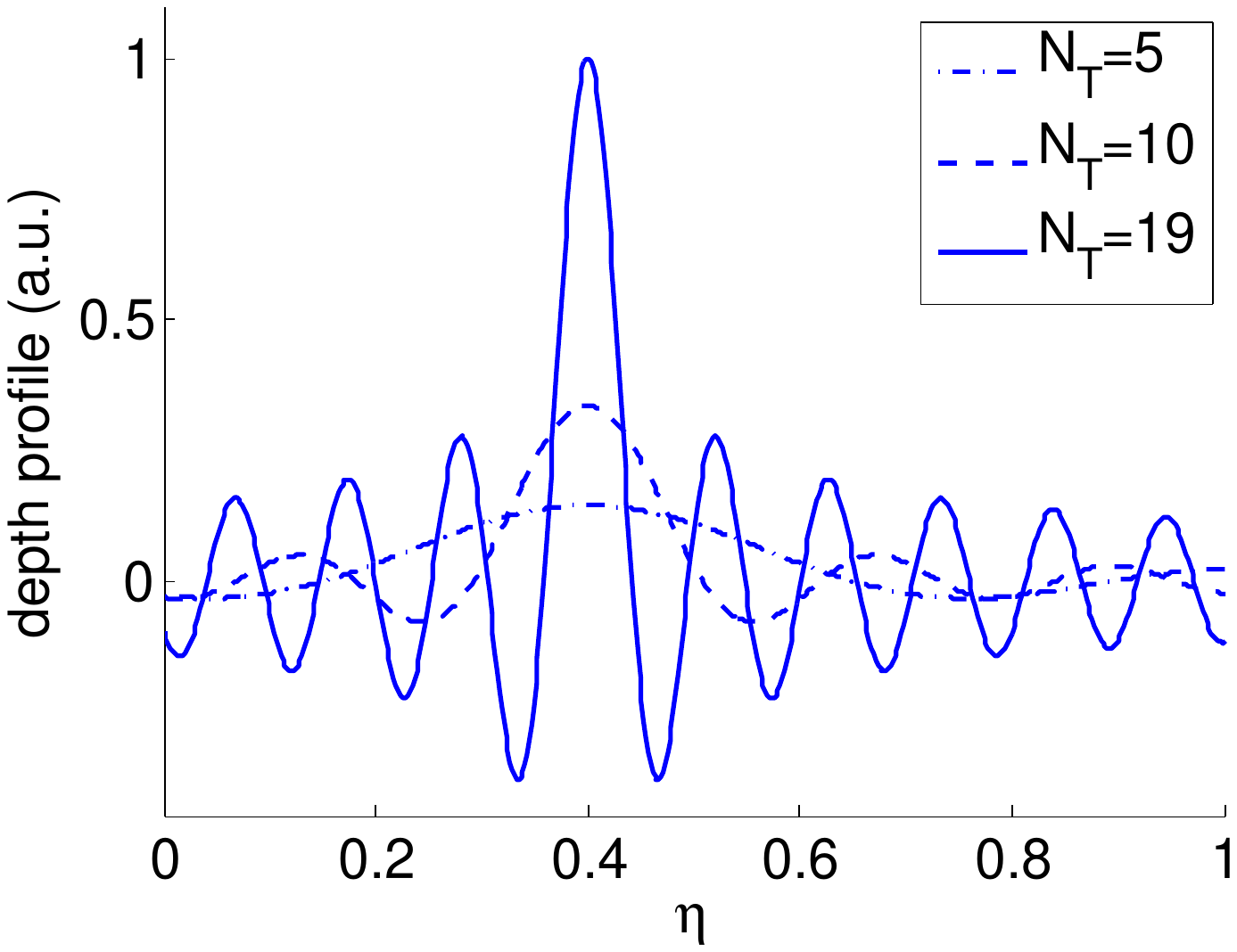}
   \includegraphics[width=8cm]{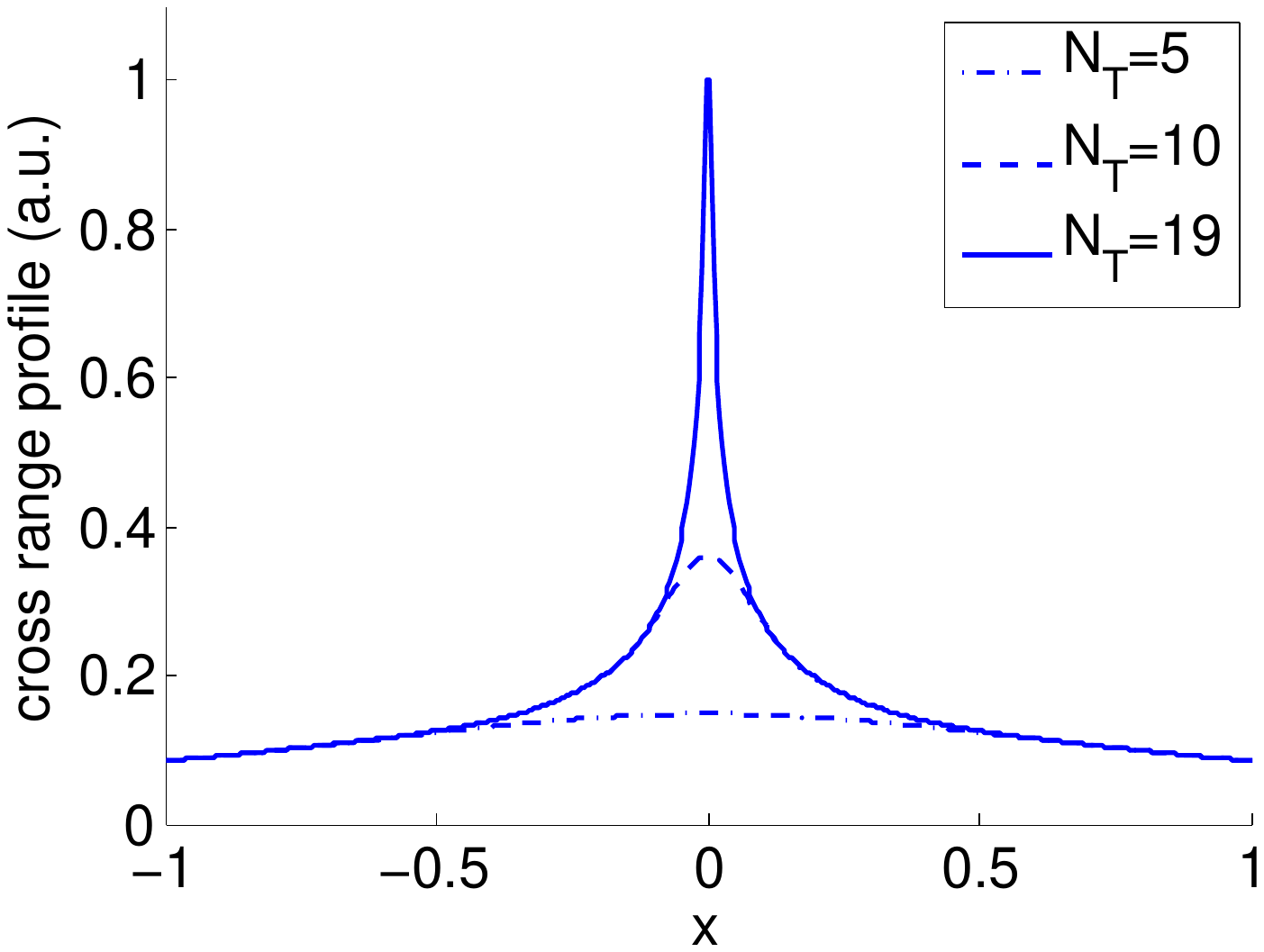}
   \end{center}
   \vspace*{-1.3in}
   \caption{Depth profile (left) and cross range profile (right) of
     the mean point spread function for the time reversal functional.
     Here $Z_\cA=100$, $\ell=1$, $\sigma=0.25$, $k=60$, ${\mathcal
       D}=1$ (so that $N=19$). The array diameter $ |\cA_X|$ is
     supposed to be smaller than the critical value (\ref{eq:TR.20})
     which is about $220$.  $N_T$ is the cut-off number (modes smaller
     than $N_T$ are recorded and reemitted). Note that the high modes
     play an important role. The larger $N_T$ is, the better the resolution.}
   \label{fig_parax1} 
\end{figure}

Consider for simplicity the full aperture case, where
\begin{eqnarray}
  \EE\left[\FTR(\om,X^\star,\eta^s)\right] &\approx& \frac{|\cA_X|}{8
  \pi Z_\cA} \sum_{j=1}^{N_T} \frac{2}{{{\mathcal D}} \beta_j} \cos
  \left[ \Big( j-\frac{1}{2}\Big)\frac{k\eta^s}{N}\right] \cos \left[
  \Big( j-\frac{1}{2}\Big)\frac{k\eta^\star}{N}\right] \nonumber \\
  &\approx& \frac{|\cA_X|}{8 \pi^2 Z_\cA N} \sum_{j=1}^{N_T}
  \frac{\cos \left[ \frac{(j-1/2)}{N} k(\eta^\star-\eta)\right]}{
  \left[ 1 - \frac{(j-1/2)^2}{N^2}\right]^{1/2}} \nonumber \\
  &\approx& \frac{|\cA_X|}{8 \pi^2 Z_\cA}\Lambda_{N_T/N}
  \big(k(\eta^s-\eta^\star) \big), \quad \quad \Lambda_\alpha(x) =
  \int_0^{\alpha} d s \, \frac{\cos ( s x )}{\sqrt{1-s^2}}.
\end{eqnarray}
The function $\Lambda_\alpha$ becomes proportional to the Bessel function
 of first kind $J_0$ as $\alpha \to 1$, more explicitly, we have
 $\Lambda_1(x) = (\pi/2) J_0(x)$ so that
\begin{eqnarray}
  \EE\left[\FTR(\om,X^\star,\eta^s)\right] \approx
  \frac{|\cA_X|}{16 \pi Z_\cA} J_0\left[ k(\eta^s-\eta^\star)\right],
  \qquad \mbox{if} ~ ~N_T \approx N.
\end{eqnarray}
We can then estimate the depth resolution as the distance between the
peak of $J_0$, that occurs when $\eta^s = \eta^\star$, and its first zero,
that occurs when $ k |\eta^s-\eta^\star| \approx 2.4$ (first zero of
$J_0$). Therefore, the depth resolution of time reversal with full
aperture is equal to the diffraction limit
\begin{equation}
|\eta^s-\eta^\star| \le \Delta_\eta(\om) = \frac{2.4}{k},
\end{equation}
if the array records the waves long enough to capture almost all the
propagating modes. The resolution deteriorates if $N_T$ is much
smaller than $N$. Indeed for small $\alpha$ we have $\Lambda_\alpha (x)
\approx (\alpha/\pi)\, {\rm sinc}(\alpha x)$ and therefore the depth
resolution is $\Delta_\eta(\om)\approx \pi N /(k N_T)$ (see Figure \ref{fig_parax1}).

\subsection{Statistical stability}
\label{sect:TRStatStab}
We now show that the time reversal function is statistically stable,
meaning that the refocusing of the wave at the original source
location does not depend on the the realization of the random medium
but only on its statistical distribution, and the point spread
function is approximately equal to its expectation.  

We restrict the analysis of statistical stability to the case of full 
aperture, where the calculations are simpler because the coupling 
matrix $\Gamma_{jl}$ becomes the identity. The point spread function 
follows from (\ref{eq:TR.12p})
\begin{eqnarray*}
 \FTR(\om,X^s,\eta^s) = \sum_{j=1}^N
\frac{\phi_j(\eta^s) \phi_j(\eta^\star)}{\beta_j(\om)} \psi
\left(\frac{ \beta'_j(\om) Z_\cA}{T}\right)
 \int dX \, \frac{1_{\cA_X}(X)}{|\cA_X|} \cT_j(\om,X,X^s,Z_\cA)
\overline{\cT_j(\om,X,X^\star,Z_\cA)},
\end{eqnarray*}
and its variance at the source location is 
\begin{eqnarray*}
{\rm Var}\left[ \FTR(\om,X^\star,\eta^\star) \right]= \sum_{j,J=1}^N
\frac{\phi_j^2(\eta^\star) }{\beta_j(\om)} \frac{\phi_J^2(\eta^\star)
}{\beta_J(\om)} \psi \left(\frac{ \beta'_j(\om)
Z_\cA}{T}\right) \psi \left(\frac{ \beta'_J(\om)
Z_\cA}{T}\right) \iint dX d Y \, \frac{1_{\cA_X}(X)1_{\cA_X}(Y)}{|\cA_X|^2}
\nonumber \\ \times \Big\{ \EE
\big[\left| \cT_j(\om,X,X^\star,Z_\cA)\right|^2 \left|
\cT_J(\om,Y,X^\star,Z_\cA)\right|^2 \big]-\EE \big[
\left|\cT_j(\om,X,X^\star,Z_\cA)\right|^2 \big] \EE\big[
\left|\cT_J(\om,Y,X^\star,Z_\cA)\right|^2 \big] \Big\} .
\end{eqnarray*}
From appendix \ref{ap:FM} (first case) we find that the variance is
much smaller than the square expectation when $|A_X|\gg \ell$, and
therefore the point spread function is equal to its mean
approximately.  The results contained in appendix \ref{ap:FM} (second
case) also show that if $|A_X| < \ell$, the variance of the point
spread function is large, and therefore that time reversal refocusing
may be unstable in this case.

There is however another mechanism that can ensure statistical
stability of the focal spot if the array is small. Indeed, if the
bandwidth of $\varphi$ is larger than the decorrelation frequency,
then the variance
\begin{eqnarray*}
{\rm Var} \left[ \ITR(X^\star,\eta^\star) \right]= \int \frac{d \om}{2
\pi}\int \frac{d \om'}{2 \pi} \overline{\hat \varphi(\om,X)} \, \hat
\varphi(\om')  \, {\rm Cov }\big[
\FTR(\om,X^\star,\eta^\star),
\FTR(\om',X^\star,\eta^\star)\big]
\end{eqnarray*}
is small because the covariance of the point spread function at two
frequencies becomes approximately zero if the frequency gap is large
enough.  Therefore, if the pulse has large bandwidth, then the
time-reversal focal spot is statistically stable even for small
arrays.

\section{Imaging}
\label{sect:IM}

The sharp and stable focusing of the time reversal process in the
random waveguide is due the backpropagation of the time reversed field
$\DTR$ in \emph{exactly the same waveguide}.  Time reversal is a
physical experiment, where the waves  can be observed in the vicinity
of the source, as they refocus. In imaging we only have access to
the data measured at the array, and the backpropagation to the search
points is synthetic. Because we cannot know the fluctuations of the
boundary, we simply ignore them in the synthetic backpropagation and
obtain the so-called \emph{reversed time migration} imaging function.
We analyze it in section \ref{sect:KM} and show that it does not give
useful results in the low-SNR regime. In particular, we show that the
images are not statistically stable with respect to realizations of
the fluctuations. Stability can be achieved by imaging with local
cross-correlations of the array measurements. Local means that we
recall the decorrelation properties of the random mode amplitudes
described in section \ref{sect:STAT.3}, and cross-correlate the
measurements over receivers located at nearby cross-ranges $X$, and
projected on the same eigenfunctions.  The resulting \emph{coherent
  interferometric} imaging method is analyzed in section
\ref{sect:CINT}.
\subsection{Reverse time migration}
\label{sect:KM}
The reverse time migration function is given by the time reversed data
$\DTR$ propagated (migrated) in the ideal waveguide to the search
points $ (x^s,\eta^s,z^s) =
\left(\frac{X^s}{\eps},\eta^s,\frac{Z^s}{\eps^2} \right).  $ Its
mathematical expression follows from (\ref{eq:HOM.parax}), with
amplitudes (\ref{eq:HOM.15}) replaced by
\begin{eqnarray}
  a_{j,o}(\om,X,Z) \leadsto \int dX' \cT_{j,o}(\om,X,X',Z)
  \frac{1}{2 i \beta_j(\om)} \int_0^{\mathcal D} d \eta \,
\phi_j(\eta) \hDTR(\om,X',\eta).
\label{eq:KM.3}
\end{eqnarray}
The ideal transfer coefficients
\begin{equation}
\cT_{j,o}(\om,X,X',Z) = \sqrt{\frac{\beta_j(\om)}{2 \pi i Z}}
  \exp \left[\frac{i \beta_j(\om) (X-X')^2}{2 Z}
\right]
\label{eq:KM.4}
\end{equation}
are defined by the Green's functions of the paraxial operator in
(\ref{eq:HOM.16}).  We obtain
\begin{eqnarray}
  \IKM(X^s,\eta^s,Z^s) = \sum_{j=1}^N \phi_j(\eta^s) \int \frac{d
  \om}{2 \pi} \, \frac{\exp \left[ i \beta_j(\om)
  \frac{Z_\cA-Z^s}{\eps^2} -i \om t \right]}{2 i \beta_j(\om)} \int d
  X \int d \eta \, \phi_j(\eta) \nonumber \\ \times \,
  \hDTR(\om,X,\eta) \, \cT_{j,o}(\om,X,X^s,
  Z_\cA-Z^s)\Big|_{t=T^\eps},
\label{eq:KM.1}
\end{eqnarray}
with the right hand side evaluated at the same time $t=T_\eps$ as in
time reversal.

We assume again a tightly supported source density normalized by
(\ref{eq:rrho}) and substitute the model (\ref{eq:TR.7}) of $\hDTR$
in (\ref{eq:KM.1}) to obtain
\begin{eqnarray}
  \IKM(X^s,\eta^s,Z^s) =
  \int \frac{d \om}{2 \pi} \, \overline{\hat \varphi(\om)} \,
\FKM\left(\om,X^s,\eta^s,Z^s\right),
 \label{eq:KM.2}
\end{eqnarray}
with frequency-dependent kernel (point spread function)
\begin{eqnarray}
  \FKM\left(\om,X^s,\eta^s,Z^s\right) \approx
  \sum_{j,l=1}^N \Gamma_{jl}
  \frac{\phi_j(\eta^s) \phi_l(\eta^\star)}{4
    \beta_j(\om) \beta_l(\om)}
  \int \frac{du}{2 \pi} \, \overline{\hat \psi(u)} \int d X \,
1_{\cA_X}(X)  \nonumber \\
  \times  \cT_{j,o} \Big(\om,X,X^s,Z_\cA-Z^s \Big)
  \overline{\cT_l\Big(\om-\frac{\eps^2 u}{T},X,X^\star,Z_\cA \Big)}
  \nonumber \\
  \times \exp \left[i \beta_j(\om) \frac{Z_\cA-Z^s}{\eps^2}
    - i \beta_l\Big( \om-\frac{\eps^2 u}{T}\Big) \frac{Z_\cA}{\eps^2}
  \right].
  \label{eq:KM.2p}
\end{eqnarray}

\subsubsection{The mean imaging function}
Let us take for simplicity the case of full aperture in depth,
where the coupling matrix $\Gamma_{jl}$ given by (\ref{eq:TR.10})
becomes the identity. We obtain from (\ref{eq:KM.2p}) and the moment
formula (\ref{eq:STAT.2}) that
\begin{eqnarray}
  \EE\left[ \FKM(\om,X^s,\eta^s,Z^s)\right] =
  \frac{|\cA_X|}{8 \pi Z_\cA} \sum_{j=1}^N
  \frac{\phi_j(\eta^s) \phi_j(\eta^\star)}{\beta_j(\om)}   \psi \left(\frac{
      \beta'_j(\om) Z_\cA}{T}\right)\exp \left[
-\frac{Z_\cA}{\cS_j(\om)} - i \beta_j(\om) \frac{Z^s}{\eps^2}\right]
    \nonumber \\
    \times  \int dX \,  \frac{1_{\cA_X}(X)}{|\cA_X|}
    \exp \left[ -\frac{i \beta_j(\om)}{Z_\cA}\Big(X-\frac{X^s + X^\star}{2}
      \Big) (X^s-X^\star)\right]. ~~
\label{eq:KM.6}
\end{eqnarray}
Moreover, assuming the aperture $\cA_X$ defined in (\ref{eq:DefAX}), and 
integrating in $X$, we get 
\begin{eqnarray}
  \EE\left[\FKM(\om,X^s,\eta^s,Z^s)\right] = \frac{|\cA_X|}{8 \pi
  Z_\cA} \sum_{j=1}^N \frac{\phi_j(\eta^s)
  \phi_j(\eta^\star)}{\beta_j(\om)} \psi \left(\frac{ 
  \beta'_j(\om) Z_\cA}{T}\right)\exp \left[ -\frac{Z_\cA}{\cS_j(\om)} -
  i \beta_j(\om) \frac{Z^s}{\eps^2}\right] \nonumber \\ \times
  \mbox{\rm sinc} \left[ \frac{\beta_j(\om)|\cA_X|}{2 Z_\cA}
  (X^s-X^\star) \right] \exp \left\{ \frac{i
  \beta_j(\om)}{2Z_\cA}\left[(X^s)^2 - (X^\star)^2 \right]\right\}. ~~
\label{eq:KM.7}
\end{eqnarray}
This result is almost the same as in the ideal waveguide, except for
the damping coefficients $\exp\left[-Z_\cA/\cS_j\right].$ 

The sinc kernel in the mean point spread function gives the focusing
in cross-range, with mode-dependent resolution
\begin{equation}
  |X^S - X^\star| \le \Delta_{X,j}(\om)=
\frac{2 \pi Z_\cA}{\beta_j(\om) |\cA_X|}.
\end{equation}
The best resolution is for the first mode, that has the largest
wavenumber $ \beta_1(\om) \approx {\pi N}/{{\mathcal D}} \approx k, $ and
gives the Rayleigh cross-range resolution
\begin{equation}
\Delta_{X,1}(\om) \approx \frac{2 \pi Z_\cA}{k |\cA_X|}.
\end{equation}

The focusing of the point spread function $\FKM$ in range can only be
due to the summation of the rapidly oscillating terms $\exp\left[-i
  \beta_j Z^s/\eps^2\right]$. But these terms are weighted by
$\mbox{exp}[- Z_\cA/\cS_j]$, which decay fast in $j$. The first 
term dominates in 
\begin{eqnarray}
  \EE\left[\FKM(\om,X^\star,\eta^\star,Z^s)\right] = \frac{|\cA_X|}{8 \pi
  Z_\cA} \sum_{j=1}^N \frac{
  \phi_j^2(\eta^\star)}{\beta_j(\om)} \psi \left(\frac{ 
  \beta'_j(\om) Z_\cA}{T}\right)\exp \left[ -\frac{Z_\cA}{\cS_j(\om)} -
  i \beta_j(\om) \frac{Z^s}{\eps^2}\right],
\label{eq:KM.7.R}
\end{eqnarray}
so the mode diversity does not lead to focusing in range, as is the
case in ideal waveguides. Nevertheless, the mean reverse time
migration function peaks at $Z^s = 0$ because of the integral
over the bandwidth in
\begin{equation}
\EE \left[ \IKM(X^\star,\eta^\star, Z^s)\right] = \int \frac{d \om}{2 \pi} 
\overline{\hat \varphi(\om)} \EE\left[\FKM(\om,X^\star,\eta^\star,Z^s)\right],
\end{equation}
and the range resolution is  of the order $\eps^2/[\beta'_1(\om_o) B]$.

When we evaluate the point spread function at 
$Z^s = 0$ and $X^s = X^\star$, we obtain
\begin{eqnarray}
  \EE \left[ \FKM(\om,X^s= X^\star,\eta^s,Z^s= 0)\right] =
  \frac{|\cA_X|}{8 \pi Z_\cA} \sum_{j=1}^N
  \frac{\phi_j(\eta^s) \phi_j(\eta^\star)}{\beta_j(\om)}   \psi \left(\frac{
      \beta'_j(\om) Z_\cA}{T}\right)
  e^{-{Z_\cA}/{\cS_j(\om)}}. \qquad
\label{eq:KM.9}
\end{eqnarray}
This is a sum of the oscillatory functions
\[
\phi_j(\eta^s) \phi_j(\eta^\star) = \frac{1}{\mathcal D} \left\{ \cos
  \left[ \pi\left(j-\frac{1}{2}\right) \frac{(\eta^s-\eta^\star)}{\mathcal
      D}\right] + \cos
  \left[ \pi\left(j-\frac{1}{2}\right) \frac{(\eta^s+\eta^\star)}{\mathcal
      D}\right]\right\}
\]
multiplied by positive weights, which are small and decay fast in $j$.
The first term dominates in (\ref{eq:KM.9}) and there is no depth
resolution at all.  We show next that these small weights also
indicate the lack of statistical stability of the reverse time
migration function.

\subsubsection{Stability analysis} To assess the stability of the 
reverse time migration, we calculate its variance at the source
location
\[
\mbox{Var}\left[ \IKM(X^\star,\eta^\star,0)\right] =
\EE\left[\left|\IKM(X^\star,\eta^\star,0)\right|^2\right] -
\left|\EE\left[\IKM(X^\star,\eta^\star,0)\right]\right|^2.
\]
We have from the results above that 
\begin{eqnarray}
  \EE\left[\IKM(X^\star,\eta^\star,0)\right] \approx
\frac{|\cA_X|}{8 \pi Z_\cA}
  \int  \frac{d \om}{2 \pi} \,
  \overline{\hat \varphi(\om)} \sum_{j=1}^N
  \frac{\phi_j^2(\eta^\star)}{\beta_j(\om)}
  \psi \left(\frac{
      \beta'_j(\om) Z_\cA}{T}\right)e^{
    -{Z_\cA}/{\cS_j(\om)}}.
\label{eq:SNRKM.2}
\end{eqnarray}
The second moment of $\IKM$ is
\begin{eqnarray}
  \EE\left[\left|\IKM(X',\eta',0)\right|^2\right] \approx
  \int \frac{d \om_1}{2 \pi} \int \frac{d \om_2}{2 \pi}
  \, \overline{\varphi(\om_1)} \varphi(\om_2)
  \sum_{j,l=1}^N \frac{\phi_j^2(\eta^\star)\phi_l^2(\eta^\star)}{16
    \beta_j^2(\om_1)\beta_l^2(\om_2)} \int \frac{d u_1}{2 \pi}
  \int \frac{d u_2}{2 \pi} \,
  \overline{\hat \psi(u_1)}  {\hat \psi(u_2)} ~~\nonumber \\
  \times \, \exp \left[ \frac{i [\beta'_j(\om_1) u_1 -
\beta'_l(\om_2)
      u_2]Z_\cA}{T} \right] \int d X_1 \int d X_2 \,
  1_{\cA_X}(X_1) 1_{\cA_X}(X_2) \, \cT_{j,o}(\om_1,X_1,X^\star,Z_\cA)~~
  \nonumber \\
  \times \,  \overline{\cT_{l,o}
    (\om_1,X_1,X^\star,Z_\cA)}  \EE \left[ \overline{\cT_{j}
      \left(\om_1-\frac{\eps^2 u_1}{T},X_1,X^\star,Z_\cA\right)} \cT_{l}
    \left(\om_2-\frac{\eps^2 u_2}{T},X_2,X^\star,Z_\cA\right)\right],~~
\end{eqnarray}
and we recall from Proposition \ref{prop.3} that only the diagonal
terms $j = l$ contribute to the expectation.  We also assume a small
bandwidth $B \ll \Omega_{d,j}$, for all the modes $j$, so that we can
use the simpler moment formula (\ref{eq:2M.4}). We obtain
\begin{eqnarray}
  \EE\left[\left|\IKM(X',\eta',0)\right|^2\right] \approx
  \frac{|\varphi(0)|^2}{(8 \pi Z_\cA)^2}
  \sum_{j=1}^N \frac{\phi_j^4(\eta^\star)}{
    \beta_j^2(\om_o)} \left|\psi \left(\frac{
        \beta'_j(\om_o) Z_\cA}{T}\right)\right|^2
  \int d X_1 \int d X_2 \,
  1_{\cA_X}(X_1) \nonumber \\
  \times  \, 1_{\cA_X}(X_2) \exp \left[ - \frac{(X_1-X_2)^2}{
2 X_{d,j}^2(\om_o)}\right],
\label{eq:KM2nd}
\end{eqnarray}
after approximating the modal wavenumbers by their value at the
central frequency. This expression can be approximated further, after
integrating in $X_1$ and $X_2$, and supposing that the decoherence
lengths $X_{d,j}$ are much smaller than the array aperture,
\begin{eqnarray}
  \EE\left[\left|\IKM(X',\eta',0)\right|^2\right] \approx
  \frac{|\varphi(0)|^2 |\cA_X| \sqrt{2 \pi}}{(8 \pi Z_\cA)^2}
  \sum_{j=1}^N \frac{X_{d,j}(\om_o)\phi_j^4(\eta^\star)}{
    \beta_j^2(\om_o)} \left|\psi \left(\frac{
        \beta'_j(\om_o) Z_\cA}{T}\right)\right|^2.
\label{eq:KM2nd.p}
\end{eqnarray}

The second moment (\ref{eq:KM2nd.p}) is clearly much larger than the
square of the mean (\ref{eq:SNRKM.2}), which is exponentially small in
range. Although the mean of the imaging function is focused at the
source, it cannot be observed because it is dominated by its random
fluctuations.  The reverse time migration lacks statistical stability
with respect to the realizations of the random fluctuations of the
boundary of the waveguide.

The calculations above are for a small bandwidth, satisfying $B \ll
\Omega_{d,j}$ for all the modes captured in the recording window.  The
calculations are more complicated for a larger bandwidth, but the
conclusion remains that reverse time migration is not stable with
respect to different realization of the random boundary fluctuations.

\subsection{Coherent interferometric imaging}
\label{sect:CINT}
The main idea of the coherent interferometric (CINT) imaging approach
is to backpropagate synthetically to the imaging points the local
cross-correlations of the array measurements, instead of the
measurements themselves. By local we mean that because of the
statistical decorrelation properties of the random mode amplitudes
described in section \ref{sect:STAT.3}, we cross-correlate the data
$\hat D(\om,X,\eta)$ at nearby frequencies and cross-ranges $X$, after
projecting it on the subspace of one eigenfunction $\phi_j$ at a time.
The projection gives the coefficients
\begin{equation}
\hat D_j(\om,X) = \int_0^{\cD} d \eta \, \phi_j(\eta) \hat D(\om,X,\eta),
\label{eq:CINT.1}
\end{equation}
which are directly proportional to the coefficients $\hat F_j$ of the
source only in the case of an array spanning the entire depth of
the waveguide. We assume this case here, because it simplifies
the analysis of the focusing and stability of the CINT function.  We
also take a small source, meaning that we essentially compute the CINT
point spread function.

The model of the coefficients (\ref{eq:CINT.1}) is
\begin{eqnarray}
  \hat D_j(\om,X) \approx 1_{\cA_X}(X) \frac{\hat
  \varphi(\om)\phi_j(\eta^\star )}{2 i \beta_j(\om)} e^{i
  \beta_j(\om)Z_\cA/\eps^2} \int \frac{d u}{2 \pi} \, \hat \psi(u)
  e^{- i \beta'_j(\om)u Z_\cA/T} \cT_j \Big(\om-\frac{\eps^2
  u}{T},X,X^\star,Z_\cA\Big),
\label{eq:CINT.2}
\end{eqnarray}
and we cross-correlate them at cross-ranges satisfying $ |X_1 - X_2|
\le X_{d,j}(\om), $ and at frequency offsets 
\begin{equation} 
|\om_1-\om_2| \le
\Omega \ll \Omega_{d,j}.  
\label{eq:Omega}
\end{equation}
We take such small $\Omega$ to simplify the second moment formulas.

The CINT image is formed by backpropagating the cross-correlations to
the imaging point, using the Green's function in the ideal waveguide.
We first define the CINT image in the $(X,Z)$-domain:
\begin{equation}
\label{eq:CINT.5a}
  \ICINT(X^s , Z^s) = \sum_{j=1}^N   \ICINT_j(X^s,Z^s)
\end{equation}
with
\begin{eqnarray}
  \ICINT_{j}(X^s, Z^s) &=&
  \iint \frac{d \om_1}{2 \pi} \frac{d \om_2}{2
    \pi} 1_\Omega(\om_1-\om_2) e^{i [\beta_j(\om_2)-\beta_j(\om_1)]     \frac{Z^s- Z_\cA}{\eps^2}} 
  \iint dX_1 d X_2 \,
  1_{X_{d,j}}(X_1-X_2)   \nonumber \\
   &&\times
  \, \hat D_j (\om_1,X_1)\overline{\hat D_j (\om_2 ,X_2 )}
  \overline{\cT_{j,o}\Big(\om_1,X_1,X^s, Z_\cA-Z^{s}\Big)}
 \cT_{j,o}\Big(\om_2,X_2,X^s, Z_\cA-Z^s\Big) , \quad
\label{eq:CINT.5}
\end{eqnarray}
%
where $1_{X_{d,j}}$ are indicator functions of the cross-range
interval $\left[-X_{d,j}(\om),X_{d,j}(\om)\right]$ calculated at the
central frequency $ \om = (\om_1+\om_2)/{2}.  $ Similarly, $1_\Omega$
is the indicator function of the frequency interval
$[-\Omega,\Omega]$.

\subsubsection{The mean CINT function}
To study the focusing of CINT, we consider its expectation
\begin{eqnarray}
  \EE\left[\ICINT(X^s,Z^s) \right] \approx \int \frac{d \om}{2
  \pi} \left| \hat \varphi(\om)\right|^2
  \EE\left[\FCINT(\om,X^s,Z^s)\right],
\label{eq:CINT.6}
\end{eqnarray}
with frequency-dependent kernel
\begin{eqnarray}
\EE\left[\FCINT(\om,X^s,Z^s)\right] \approx \sum_{j=1}^N
\frac{\phi_j^2(\eta^\star)}{32 \pi^3 Z_\cA
(Z_\cA-Z^s)} \left|\psi \left(\frac{ \beta'_j(\om)
Z_\cA}{T}\right)\right|^2
\hspace{-0.05in}\iint dX_1 d X_2 \, 1_{\cA_X} (X_1)1_{\cA_X}(X_2)
\nonumber \\ \times \exp \left\{ i
\beta_j(\om)\left[\frac{(X_1-X^\star)^2}{2
Z_\cA}-\frac{(X_1-X^s)^2}{2(Z_\cA - Z^s)} - \frac{(X_2-X^\star)^2}{2
Z_\cA}+\frac{(X_2-X^s)^2}{2(Z_\cA - Z^s)}\right] -
\frac{(X_1-X_2)^2}{2 X_{d,j}^2(\om)} \right\} \nonumber \\ \times \int
\frac{d \tom}{2 \pi} \, 1_{\Omega}(\tom) \exp\left\{i
\left[\beta_j\left(\om+\frac{\tom}{2}\right)-\beta_j
\left(\om-\frac{\tom}{2}\right)\right] \frac{Z^s}{\eps^2}\right\}.~~~~
\label{eq:CINT.8}
\end{eqnarray}
This expression follows from (\ref{eq:CINT.5}), the second moment
formula (\ref{eq:2M.4}), and definition (\ref{eq:KM.4}) of the ideal
transfer coefficients $\cT_{j,o}$.

\subsubsection{Cross-range focusing}
\label{sect:CINT_CC} 
Let us consider in (\ref{eq:CINT.8}) a search point at the range of
the source $Z^s = 0$,
\begin{eqnarray}
\EE\left[\FCINT(\om,X^s,0)\right] \approx \frac{\Omega}{64 \pi^4
  Z^2_\cA} \sum_{j=1}^N \phi_j^2(\eta^\star)\left|\psi \left(\frac{
  \beta'_j(\om) Z_\cA}{T}\right)\right|^2
\hspace{-0.05in}\iint dX_1 d X_2 \, 1_{\cA_X} (X_1)1_{\cA_X}(X_2)
\nonumber \\ \times \exp \left[ i \beta_j(\om)
\frac{(X^s-X^\star)}{Z_\cA} (X_1-X_2) - \frac{(X_1-X_2)^2}{2
X_{d,j}^2(\om)} \right]. \quad 
\label{eq:CINT.9}
\end{eqnarray}
This formula simplifies after integrating over the array aperture and
assuming as before that $X_{d,j} \ll |\cA_X|$, 
\begin{eqnarray}
\EE\left[\FCINT(\om,X^s,0)\right] \sim \frac{\Omega |\cA_X| (2
  \pi)^{1/2}}{64 \pi^4 Z^2_\cA} \sum_{j=1}^N \phi_j^2(\eta^\star)
X_{d,j}(\om) \left|\psi \left(\frac{ \beta'_j(\om)
  Z_\cA}{T}\right)\right|^2 \nonumber \\ \times \, \exp \left\{ -
\frac{1}{2}\left[ \frac{\beta_j(\om) (X^s-X^\star)
    X_{d,j}(\om)}{Z_\cA} \right]^2 \right\}.
\label{eq:CINT.10}
\end{eqnarray}
Each term in the sum focuses at the source, with resolution 
\begin{equation}
|X^s - X^\star| \le \Delta_{X,j}(\om) = \frac{2 Z_\cA}{\beta_j(\om)
 X_{d,j}(\om)}
\label{eq:CINT.11}
\end{equation}
defined as twice the standard deviation of the Gaussian in
(\ref{eq:CINT.10}). The number of modes participating in the sum is
determined by the length of the recording time window, as before, but
each mode is weighted by the correlation length $X_{d,j}$, which
decreases monotonically with $j$. The first mode has the largest
contribution in (\ref{eq:CINT.10}), and gives the best cross-range
resolution. Since its wavenumber is approximately $\beta_1(\om)
\approx \pi N/\cD \approx k$,
\begin{equation}
\label{eq:CINT.12}
\Delta_{X,1}(\om) \approx \frac{ 2 Z_\cA}{k X_{d,j}(\om)} \sim 
\frac{2 \pi Z_\cA}{k \left[\pi X_{d,1}(\om)\right]}
\end{equation}
is comparable to the classic Rayleigh resolution for an array of 
aperture equal to $\pi X_{d,1}(\om)$ (see Figure \ref{fig_parax3}). 

The cross-range resolution (\ref{eq:CINT.12}) is worse than that of
time reversal. Scattering at the random boundary is beneficial to the
time reversal process, and the more modes are recorded, the better the
result. However, scattering impedes imaging, and the best cross-range
resolution is achieved with the first mode. Even with this mode, the
resolution is worse than that in ideal waveguides $2 \pi Z_\cA/(k
|\cA_X|)$, because $X_{d,1} \ll |\cA_X|$.

\subsubsection{Range focusing}
When we evaluate the mean CINT point spread function (\ref{eq:CINT.8})
at the cross-range $X^s = X^\star$, we obtain
\begin{eqnarray}
\EE\left[\FCINT(\om,X^\star,Z^s)\right] \approx
\sum_{j=1}^N \frac{\phi_j^2(\eta^\star)}{32 \pi^3
Z_\cA (Z_\cA-Z^s)} \left|\psi \left(\frac{ \beta'_j(\om)
Z_\cA}{T}\right)\right|^2
\hspace{-0.05in}\iint dX_1 d X_2 \, 1_{\cA_X} (X_1)1_{\cA_X}(X_2)
\nonumber \\ \times \exp \left\{ - i \beta_j(\om)
\frac{(X_1-X_2)Z^s}{Z_\cA(Z_\cA-Z^s)}\left(\frac{X_1+X_2}{2} -
X^\star\right) - \frac{(X_1-X_2)^2}{2 X_{d,j}^2(\om)} \right\}
\nonumber \\ \times \int \frac{d \tom}{2 \pi} \, 1_{\Omega}(\tom)
\exp\left\{-i \left[\beta_j\left(\om+\frac{\tom}{2}\right)-\beta_j
\left(\om-\frac{\tom}{2}\right)\right] \frac{Z^s}{\eps^2}\right\}.~~~~
\label{eq:CINT.13}
\end{eqnarray}
Because we integrate over $\tom$ the rapidly oscillating integrand, 
at scale $\eps^2$, we have from the method of stationary phase that 
(\ref{eq:CINT.13}) is large for 
\[
Z^s = \eps^2 \zeta^s 
\]
with $\zeta^s$ independent of $\eps$.  Recall the assumption
(\ref{eq:Omega}) of the frequency offsets.

\begin{figure}
   \vspace*{-1.2in}
      \begin{center}
   \includegraphics[width=8cm]{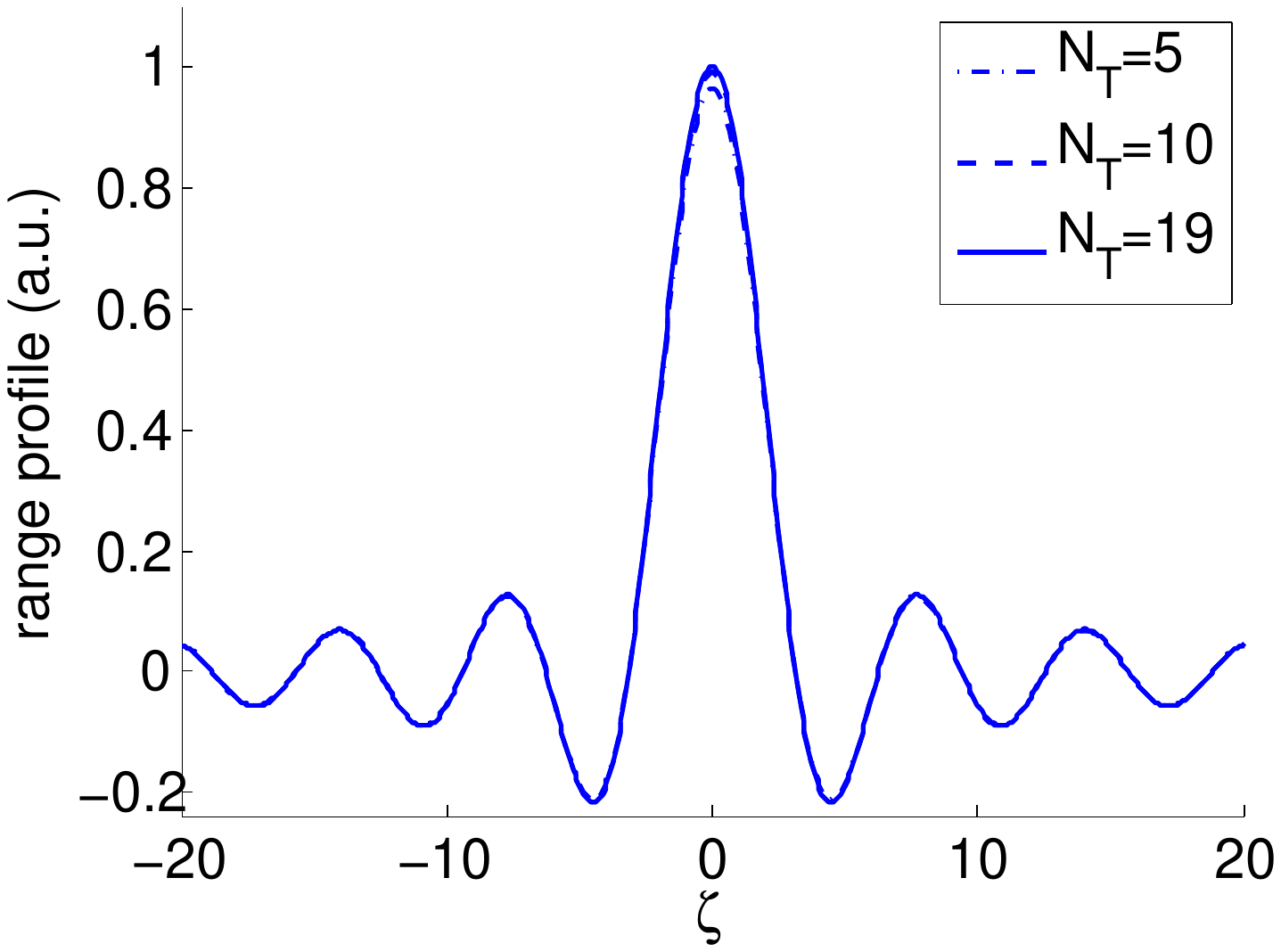}
   \includegraphics[width=8cm]{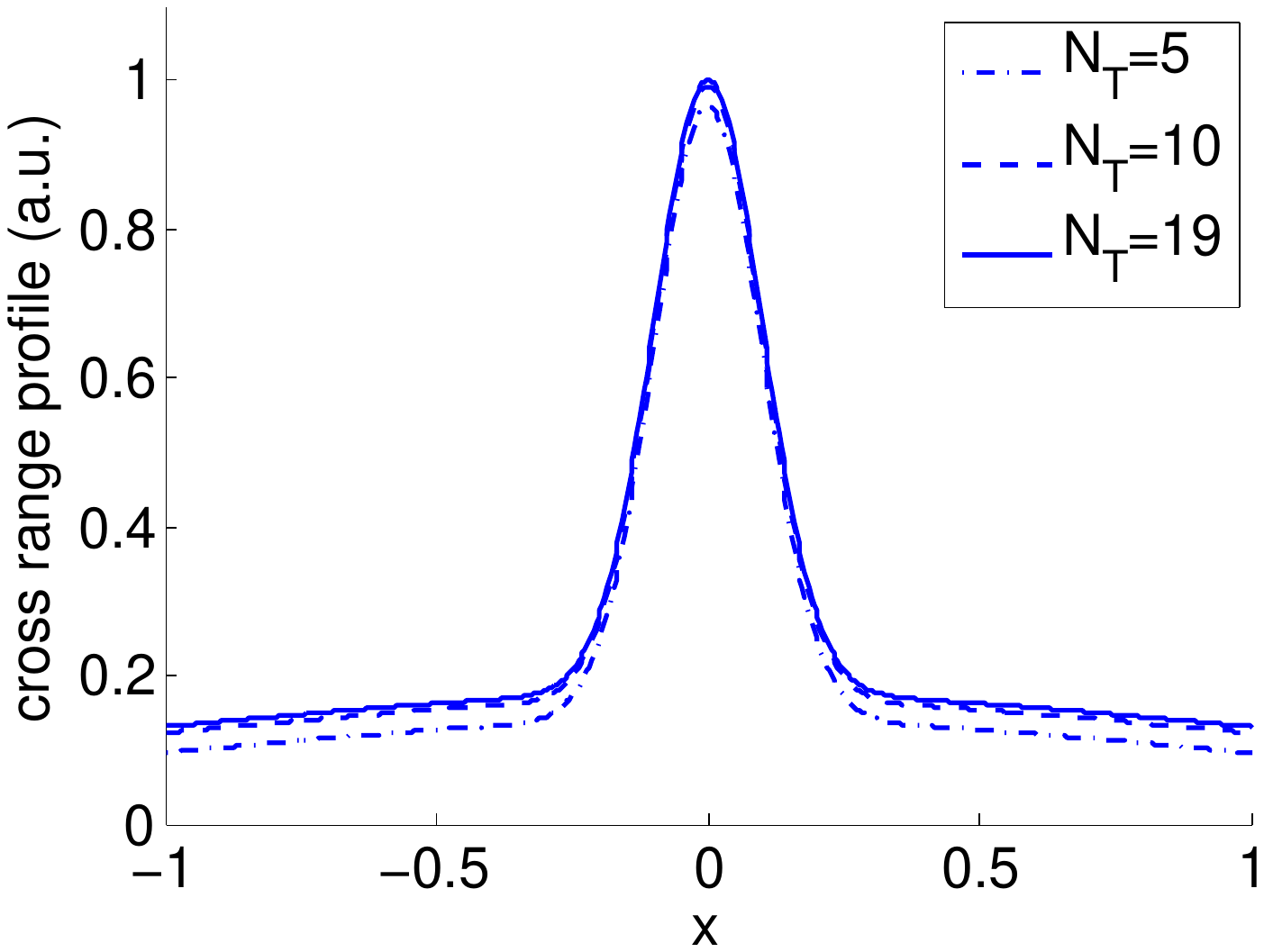}
   \end{center}
   \vspace*{-1.3in}
   \caption{Range profile and cross range profile of the mean point
     spread function for the CINT functional.  Here $Z_\cA=100$,
     $\ell=1$, $\sigma=0.25$, $k=60$, ${\mathcal D}=1$ (so that
     $N=19$), and the cut-off frequency is $\Omega/c=1$.  $N_T$ is the
     cut-off number (modes smaller than $N_T$ are recorded and
     reemitted). Note that the high modes do not play any role.}
   \label{fig_parax3} 
\end{figure}

The mean point spread function becomes 
\begin{eqnarray}
\EE\left[\FCINT(\om,X^\star, \eps^2 \zeta^s)\right] \approx
\frac{\Omega |\cA_X| (2 \pi)^{1/2}}{64 \pi^4 Z_\cA^2}\sum_{j=1}^N
X_{d,j}(\om)\phi_j^2(\eta^\star) \left|\psi \left(\frac{ \beta'_j(\om)
Z_\cA}{T}\right)\right|^2 \mbox{sinc} \left[ \beta'_j(\om) \Omega
\zeta^s\right], \quad 
\label{eq:CINT.14}
\end{eqnarray}
and we define the mode-dependent scaled range resolution by 
\begin{equation}
|\zeta^s| \le \Delta_{\zeta,j} = \frac{1}{\Omega \beta_j'(\om)}.
\end{equation}
Again, the resolution is best for the first mode, which has the
largest weight $X_{d,1}(\om)$ in (\ref{eq:CINT.14}).  See Figure
\ref{fig_parax3} for an illustration.

\subsubsection{Depth estimation}
\begin{figure}
   \vspace*{-1.2in}
      \begin{center}
   \includegraphics[width=8cm]{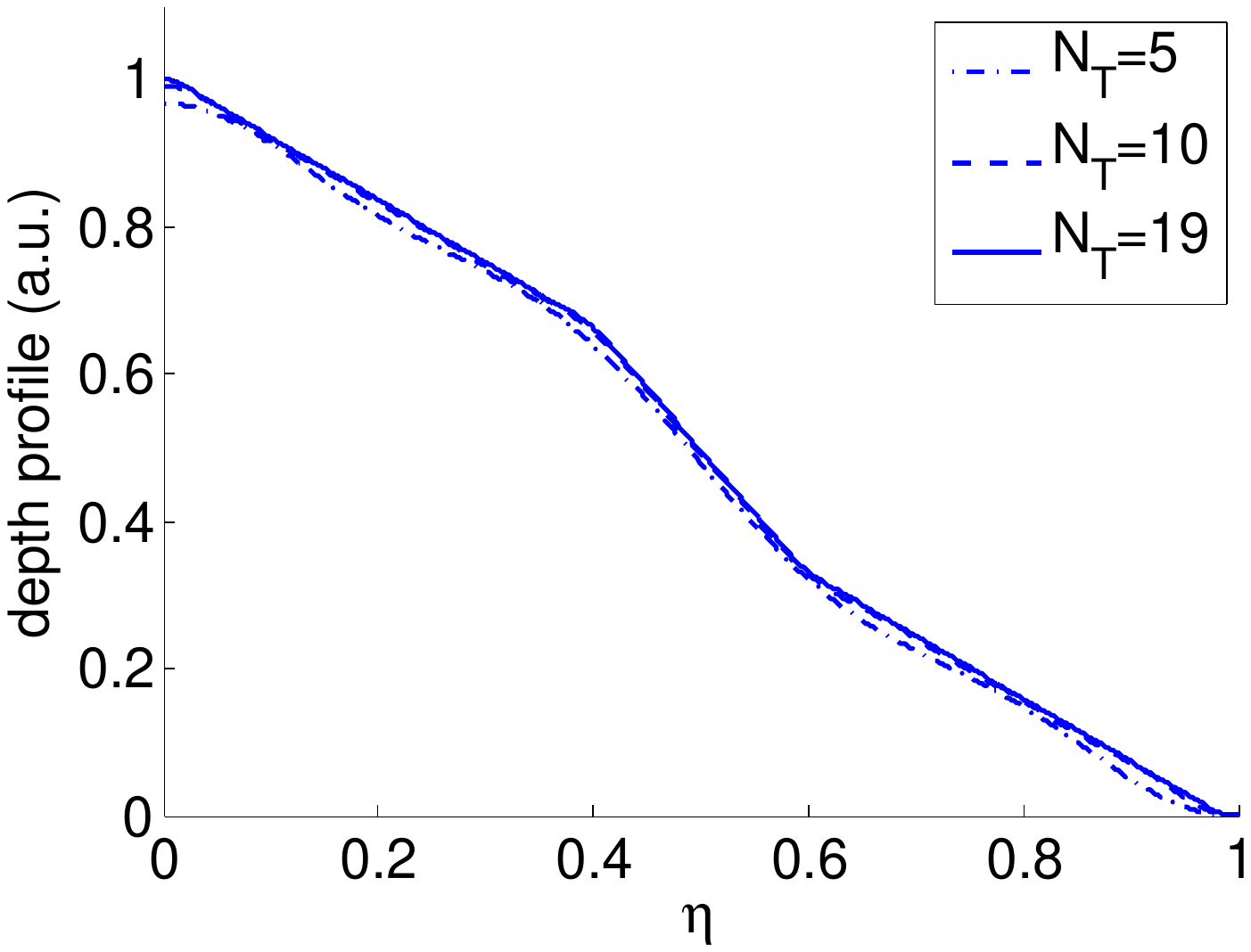}
   \includegraphics[width=8cm]{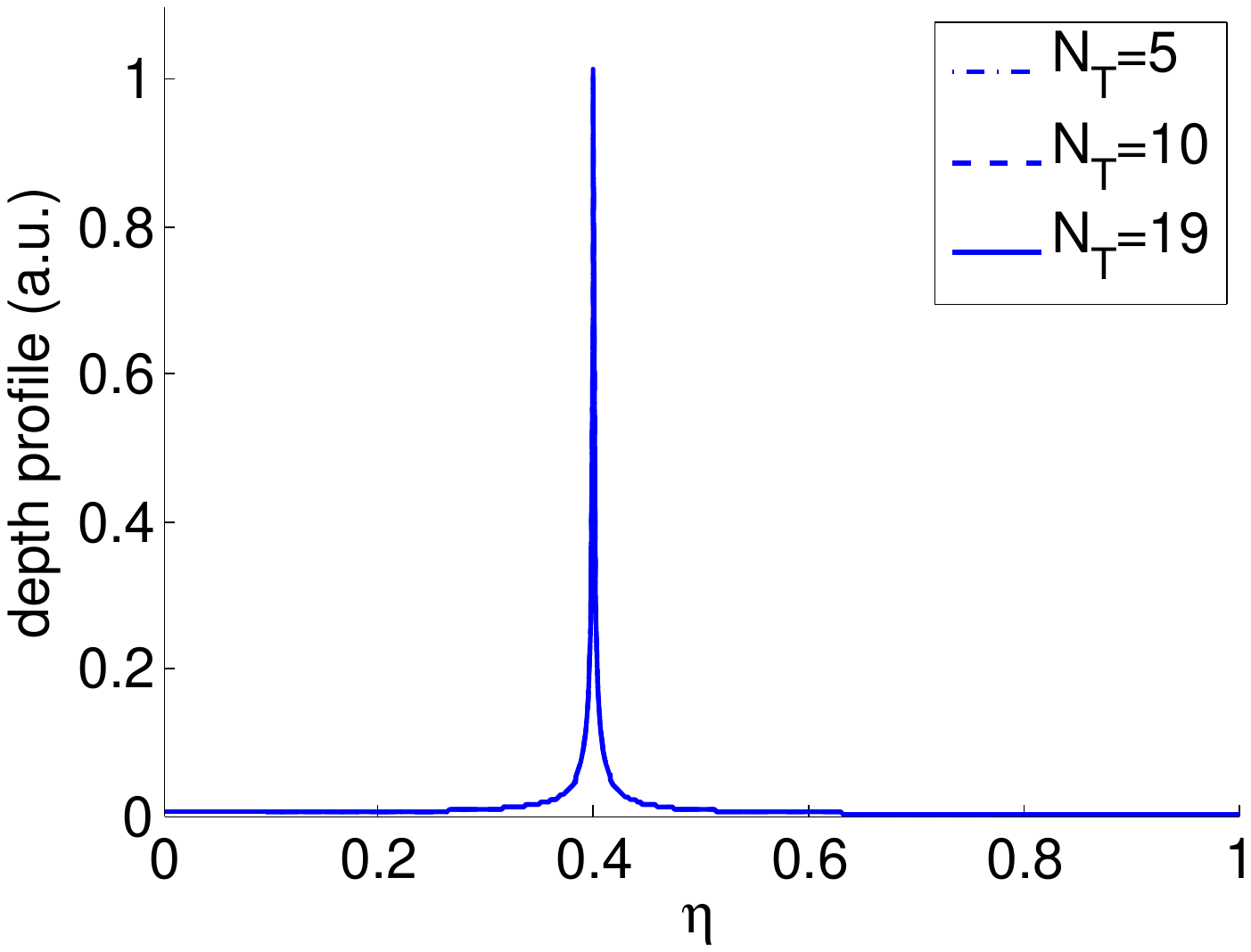}
   \end{center}
   \vspace*{-1.3in}
   \caption{Depth profile with (\ref{eq:CINT.15a}) (left) and with
     (\ref{eq:CINT.16}) (right) for the CINT functional.  In the right
     picture we plot the reciprocal of the square root of the function
     in (\ref{eq:CINT.16}).  Here $Z_\cA=100$, $\ell=1$,
     $\sigma=0.25$, $k=60$, ${\mathcal D}=1$ (so that $N=19$).  $N_T$
     is the cut-off number (modes smaller than $N_T$ are recorded and
     reemitted). Note that the high modes do not play any role.}
\label{fig_parax5}
\end{figure}

One natural way to estimate the depth $\eta^\star$ would be to consider the 
full CINT imaging functional
$$
\WICINT(X^s ,\eta^s, Z^s) = \sum_{j=1}^N   \ICINT_j(X^s,Z^s) \phi_j^2(\eta^s),
$$
with $\ICINT_j(X^s,Z^s)$ defined by (\ref{eq:CINT.5}).
However, if we define $\WFCINT$ as
\begin{eqnarray}
  \EE\left[ \WICINT(X^s,\eta^s,Z^s) \right] \approx \int \frac{d \om}{2
  \pi} \left| \hat \varphi(\om)\right|^2
  \EE\left[\WFCINT(\om,X^s,\eta^s,Z^s)\right],
\label{eq:CINT.15a}
\end{eqnarray}
and if  we take $X^s = X^\star$ and $Z^s = 0$, then 
we obtain
\begin{eqnarray}
\EE\left[ \WFCINT(\om,X^\star,\eta^s,0)\right] \approx \frac{\Omega
  |\cA_X|(2 \pi)^{1/2}}{64 \pi^4 Z_\cA^2}\sum_{j=1}^N X_{d,j}(\om)
\phi_j^2(\eta^\star) \phi_j^2(\eta^s) \left|\psi \left(\frac{
  \beta'_j(\om) Z_\cA}{T}\right)\right|^2.
\label{eq:CINT.15}
\end{eqnarray}
This is a sum of positive terms and it does not have a peak at the
depth of the source  (see Figure \ref{fig_parax5}). 

Because of scattering at the random boundary the modes are decoupled,
and we cannot speak of coherent imaging in depth. We work instead with
the squares of the mode amplitudes, i.e. intensities. Incoherent
imaging means estimating the depth of the source based on the
mathematical model (\ref{eq:CINT.15}). More explicitly, we can
estimate $\eta^\star$ by solving the least squares minimization
problem
\begin{equation}
\min_{\eta^s} \sum_{j=1}^N  \left| \ICINT_j (\hat{X}^\star,\hat{Z}^\star) - \frac{\Omega
|\cA_X|(2 \pi)^{1/2}}{64 \pi^4 Z_\cA^2}\int \frac{d\om}{2 \pi} |\hat
\varphi(\om)|^2  X_{d,j}(\om) 
\phi_j^2(\eta^s) \left|\psi \left(\frac{ \beta'_j(\om)
Z_\cA}{T}\right)\right|^2 \right|^2,
\label{eq:CINT.16}
\end{equation}
where the estimators $\hat{X}^\star$ and $\hat{Z}^\star$
of the cross-range $X^\star$ and range offset
$Z^\star = 0$ of the source
have been  determined as the location of the maximum of (\ref{eq:CINT.5a})  (see Figure \ref{fig_parax5}).

\subsubsection{Statistical stability} The analysis of statistical
stability of the CINT function is basically the same as that of time
reversal. The function is stable when evaluated in the vicinity of the
source location if the array has large aperture $|\cA_X| \gg \ell$. We
have seen in section \ref{sect:CINT_CC} that a large aperture does not
improve the focusing of $\EE\left[\ICINT\right]$. The cross-range
resolution is limited by the decoherence length. But a large aperture
is needed for the CINT function to be statistically stable.

Another way of achieving statistical stability of CINT is to have a
pulse with large bandwidth. This was already noted in the discussion
of statistical stability of time reversal in section
\ref{sect:TRStatStab}.

Note that the statistical stability of CINT relies on computing
correctly the local cross-correlations of the measurements at the
array. By this we mean that the cross-range and frequency offsets in
the correlations should not exceed the decoherence length and
frequency. Moreover, the cross-correlations should be with one mode at
a time.  This can be done with arrays that span the whole depth of the
waveguide, because the coupling matrix $\Gamma_{jl}$ becomes the
identity when $|\cA_\eta|=\cD$.  If the aperture $|\cA_\eta|$ is
small, there are large mode index offsets $|j-l|$ for which
$\Gamma_{jl} \ne 0$.  Consequently, there are many terms of the form
$\overline{\cT_j} \cT_l$, with $j \ne l$, that participate in the
expression of the imaging function. Since only the diagonal terms are
correlated, we obtain that $\ICINT$ has large variance when
$|\cA_\eta| \ll \cD$.

In practice, the decoherence scales $X_{d,j}$ and $\Omega_{d,j}$ are
likely not known explicitly. The formulas derived above are specific
to our mathematical model.  However, the decoherence scales can be
estimated as we form the image, using an adaptive procedure similar to
that introduced in \cite{ADA}.

\section{Summary}
\label{sect:SUM}
In this paper we analyze propagation of acoustic waves
in three-dimensional random waveguides.  The waves are trapped by top and
bottom boundaries, but the medium is unbounded in the remaining two
directions. The top boundary has small, random fluctuations.  We
consider a source that emits a beam, and study the resulting random
wave field in the waveguide.

The analysis is in a long range, paraxial scaling regime modeled with
a small parameter $\eps$. It is defined as the ratio of the central
wavelength $\la_o$ of the pulse emitted from the source and the
emitted beam width $r_0$.  The range of propagation is of the order of
the Rayleigh length $r_0^2/\la_o = \eps^{-2}\la_o$. The fluctuations
of the boundary are on a length scale that is similar to the
beam width, and their small amplitude is scaled so that they cause
significant cumulative scattering effects when the waves travel at
ranges of the order of the Rayleigh length.

The wave field is given by a superposition of waveguide modes with
random amplitudes. The modes are solutions of the wave equation in the
ideal waveguide, with flat boundary. The scattering effects are
captured by their random amplitudes. We show that in our scaling
regime the amplitudes satisfy a system of paraxial equations driven by
the same Brownian motion field. We use the system to calculate three
important mode-dependent scales that quantify the net scattering
effects in the waveguide, and play a key role in applications such as
imaging and time reversal. The first mode-dependent scale is the
scattering mean free path. It gives the range over which the mode
loses its coherence, meaning that the expectation of its random
amplitude is smaller than its fluctuations. The other mode-dependent
scales are the decoherence length and frequency. They give the
cross-range scale and frequency offsets over which the mode amplitudes
become statistically uncorrelated.

We use the results of the analysis to study time reversal and imaging
of the source with a remote array of sensors, in a low SNR regime. Low
SNR means that the waves travel over distances that exceed the
scattering mean free paths of all the modes, so that the random wave
field measured at the array is dominated by its fluctuations.

In time reversal, the waves received at the array are time reversed
and then re-emitted in the medium. They travel back to the source and
refocus. The refocusing is expected by the time reversibility of the
wave equation, but the resolution is limited in ideal waveguides by
the aperture of the array. We analyze the time reversal process in the
random waveguide and show that super-resolution occurs, meaning that
scattering at the random boundary improves the refocusing resolution.
An essential part of the resolution analysis is the assessment of
statistical stability with respect to different realizations of the
random boundary fluctuations.  We show that statistical stability
holds if the array has large aperture and/or the emitted pulse from
the source has a large bandwidth.

Time reversal is very different from imaging. In time reversal the
array measurements are backpropagated physically, in the real
waveguide.  In imaging we can only backpropagate the time reversed
data in software, in a surrogate waveguide.  Because we cannot know
the boundary fluctuations, we neglect them altogether, and the
surrogate is the ideal waveguide.  The resulting imaging function is
called reverse time migration and it does not work in low SNR regimes.
It lacks statistical stability, i.e., the images change unpredictably
from one realization of the fluctuations to another.

We show that robust imaging can be carried out in low SNR regimes if
we backpropagate local cross-correlations of the array measurements,
instead of the measurements themselves. Here local means that we
cross-correlate the data projected on one mode at a time, and for
nearby cross-ranges and frequencies. The method is called coherent
interferometric (CINT), because it is an extension of the CINT
approach introduced and analyzed in \cite{BPT-05,ADA,BPT-07,BGPT-11}
for imaging in open, random environments. We show that CINT images are
statistically stable under two conditions: The first condition is the
same as in time reversal and it says that the array should have a
large aperture and/or the pulse bandwidth should be large. The second
condition is that the cross-range and frequency offsets used in the
calculation of the local cross-correlations do not exceed
the mode-dependent decoherence length and frequency, respectively.  We derived
mathematical expressions of these scales, for our model. In practice,
they can be estimated adaptively, using the image formation, with an
approach similar to that in \cite{ADA}. The estimation is possible
because there is a trade-off between stability and resolution that is
quantified by the decoherence scales. If we over-estimate them we lose
statistical stability. If we under-estimate them, we lose resolution.

While cumulative scattering aids in time reversal, it impedes imaging.
We quantify this explicitly in the resolution analysis of CINT. In
time reversal the resolution improves when we record the wave field
over a long time, so that we include the high-order modes that travel
at slower speed. In CINT, the best cross-range and range resolution is
given by the first mode, which encounters the random boundary less
often, and is thus less affected by the fluctuations. The cross-range
resolution is similar to the classic Rayleigh one of range times
wavelength divided by the aperture, but instead of the real aperture
we have the decoherence length of the mode. This length decreases
monotonically with range, because longer distances of propagation in
the random waveguide mean stronger scattering effects. Similarly, the
range resolution is similar to the classic one, of speed divided by
the bandwidth, but the bandwidth is replaced by the decoherence
frequency which decreases monotonically with range. 

The estimation of the depth of the source is different than that of
range and cross-range.  Because the modes decorrelate in the low SNR
regime, we cross-correlate the data projected on one mode at a time,
so essentially, we work with intensities. The estimation of the depth
of the source from the intensities can be done by minimizing the
misfit between the processed measurements and the mathematical model.
While the cross-range and range estimation with CINT is done best with
the first waveguide mode, the depth estimation requires many modes.
Thus, we still need a long recording time at the array to capture the
later arrival of the high-order modes.

\section{Acknowledgements}
The work of L. Borcea was partially supported by the AFSOR Grant
FA9550-12-1-0117, the ONR Grants N00014-12-1-0256, N00014-09-1-0290
and N00014-05-1-0699, and by the NSF Grants DMS-0907746, DMS-0934594.
The work of J.Garnier was supported in part by ERC Advanced Grant
Project MULTIMOD-267184.

\appendix
\section{Second moment calculation}
\label{ap:SM}

\vspace{0.05in} The equation for $\EE\left[\cT_j
  (\om_1,X_1,X_1',Z)\overline{\cT_l}(\om_2,X_2,X_2',Z)\right]$ follows
from (\ref{eq:LR.14}), using It\^{o} calculus,
\begin{eqnarray}
  \partial_Z \EE\left[\cT_j \overline{ \cT_l}\right] = \left[
      \frac{i}{2 \beta_j(\om_1)} \partial_{X_1}^2 -\frac{i}{2
      \beta_l(\om_2) } \partial_{X_2}^2 - \Big( \frac{1}{
      \sqrt{\cS_j(\om_1)}} - \frac{1}{ \sqrt{\cS_l(\om_2)}}\Big)^2 -
      \frac{2C_o(X_1-X_2)}{\sqrt{\cS_j(\om_1) \cS_l(\om_2)}}
      \right]\EE\left[\cT_j \overline{\cT_l}\right]. ~ ~
\label{eq:A.01}
\end{eqnarray}
Its solution can be written as
\begin{equation}
  \EE\left[\cT_j(\om_1,X_1,X_1',Z) \overline{\cT_l(\om_2,X_2,X_2',Z)}\right] = 
  M_{jl}(\om_1,\om_2,X_1,X_2,Z; X_1',X_2') e^{- 
    \Big( \frac{1}{  \sqrt{\cS_j(\om_1)}} - 
     \frac{1}{  \sqrt{\cS_l(\om_2)}}\Big)^2\hspace{-0.05in} Z},
\label{eq:A.1}
\end{equation}
with $M_{jl}$ solving 
\begin{eqnarray}
\partial_Z M_{jl} = \left[ 
    \frac{i}{2 \beta_j(\om_1)} \partial_{X_1}^2 -\frac{i}{2 \beta_l(\om_2) }
    \partial_{X_2}^2  -
    \frac{2C_o(X_1-X_2)}{\sqrt{\cS_j(\om_1) 
        \cS_l(\om_2)}} \right]M_{jl}
\label{eq:A.2}
\end{eqnarray}
for $Z>0$, and the initial condition
\begin{equation}
M_{jl}(\om_1,\om_2,X_1,X_2,0;X_1',X_2') =\delta(X_1-X_1') \delta(X_2-X_2').
\label{eq:A.3}
\end{equation}

\subsection{Single frequency}
Let us begin with the single frequency case, $\om_1 = \om_2 = \om$,
and introduce the center and difference coordinates $\xi$ and $\txi$ so that
\begin{equation}
X_1 = \frac{\xi + \txi/2}{\sqrt{\beta_j(\om)}}, \quad 
X_2 = \frac{\xi - \txi/2}{\sqrt{\beta_l(\om)}}.
\label{eq:A.4}
\end{equation}
In this coordinate system we have that 
\begin{equation}
  U_{jl}(\om,\xi,\txi,Z;\xi',\txi') = M_{jl}\Big(\om,\om,
  \frac{\xi + \txi/2}{\sqrt{\beta_j(\om)}},\frac{\xi - \txi/2}{
    \sqrt{\beta_l(\om)}},Z; \frac{\xi' + \txi'/2}{\sqrt{\beta_j(\om)}},
\frac{\xi' - \txi'/2}{\sqrt{\beta_l(\om)}}\Big)
\label{eq:A.5}
\end{equation}
satisfies the initial value problem 
\begin{eqnarray}
  \partial_Z U_{jl} &=& i \partial_{\xi} \partial_{\txi} U_{jl} 
  - \frac{2}{\sqrt{\cS_j  \cS_l}} C_o\Big[\Big(\frac{1}{\sqrt{\beta_j}}-
  \frac{1}{\sqrt{\beta_l}}\Big) \xi + 
  \Big(\frac{1}{\sqrt{\beta_j}}+\frac{1}{\sqrt{\beta_l}}\Big) 
  \frac{\txi}{2}\Big] U_{jl}, \qquad Z > 0, \nonumber \\
  U_{jl} &=& \sqrt{\beta_j \beta_l} \, \delta(\xi-\xi') 
\delta(\txi-\txi'), \qquad Z = 0.
\label{eq:A.6}
\end{eqnarray}
Its Fourier transform in $\txi$ is the Wigner distribution 
\begin{equation}
  W_{jl}(\om,\xi,\tk,Z;\xi',\txi') = 
\int_{-\infty}^\infty \frac{d \txi}{2 \pi} \, 
  U_{jl} (\om,\xi,\txi,Z;\xi',\txi')e^{-i \tk \txi},
\label{eq:A.7}
\end{equation} 
the solution of the transport equation
\begin{eqnarray}
  \left[\partial_Z + \tk \partial_\xi\right] W_{jl}(\om,\xi,\tk,Z;\xi',\txi') =
  - \frac{4 \sqrt{\beta_l}}{\sqrt{\cS_j \cS_l} (\sqrt{\beta_j} + 
    \sqrt{\beta_l})} \int d q \, \widehat{C}_o \Big( q \frac{2 \sqrt{\beta_l}}{
    (\sqrt{\beta_j} + 
    \sqrt{\beta_l})}\Big)\nonumber \\
  \times   \exp \left[ - \frac{i q \xi}{\sqrt{\beta_j}} 
    \frac{2 (\sqrt{\beta_j}-\sqrt{\beta_l})}{(\sqrt{\beta_j}+\sqrt{\beta_l})}
  \right]
  W_{jl}\Big(\om,\xi,\tk-\frac{q}{\sqrt{\beta_j}},Z;\xi',\txi'\Big),
\label{eq:A.8}
\end{eqnarray}
for $Z > 0$, with initial condition 
\begin{equation}
  W_{jl}(\om,\xi,\tk,0;\xi',\txi') = \frac{\sqrt{\beta_j \beta_l}}{2 \pi}
 e^{-i \tk \txi'}\delta(\xi-\xi'),
\label{eq:A.9}
\end{equation}
and kernel
\begin{equation}
\widehat{C}_o(\kappa) = \delta(\kappa) - \frac{\widehat R_o(\kappa)}{R_o(0)}.
\label{eq:A.10}
\end{equation}
Here 
$$
\widehat R_o(\kappa) = \frac{1}{2\pi} \int R_o(\xi) e^{-i \kappa\cdot \xi} d\xi.
$$

\subsubsection{Single mode moments}
The transport equation (\ref{eq:A.8}) simplifies in the case $j = l$, 
\begin{eqnarray}
  \left[\partial_Z + \tk \partial_\xi\right] W_{jj}(\om,\xi,\tk,Z;\xi',\txi') 
  &=&
  - \frac{2}{\cS_j}  \int d q \, \widehat{C}_o ( q) 
  W_{jl}\Big(\om,\xi,\tk-\frac{q}{\sqrt{\beta_j}},Z;\xi',\txi'
  \Big), \nonumber \\
  W_{jj}(\om,\xi,\tk,0;\xi',\txi') &=& \frac{\beta_j} {2 \pi} e^{-i \tk \txi'}
\delta(\xi-\xi'),
\label{eq:A.11}
\end{eqnarray}
and can be integrated easily after Fourier transforming in $\tk$ and $\xi$.
Explicitly,
\begin{equation}
  V_{jj}(\om,\kappa,\txi,Z;\xi',\txi') = 
\int_{-\infty}^\infty \frac{d \xi}{2 \pi}
  \int_{-\infty}^\infty d \tk \, W_{jj}(\om,\xi,\tk,Z;\xi',\txi') 
e^{-i \kappa \xi +     i \tk \txi}
\label{eq:A.12}
\end{equation}
satisfies the initial value problem
\begin{eqnarray}
  \left[ \partial_Z + \kappa \partial_{\txi} \right]
V_{jj}(\om,\kappa,\txi,Z;\xi',\txi') &=& -\frac{2}{\cS_j}{C}_o
\Big(\frac{\txi}{\sqrt{\beta_j}} \Big)
V_{jj}(\om,\kappa,\txi,Z;\xi',\txi'), \qquad Z > 0, \nonumber \\
V_{jj}(\om,\kappa,\txi,0;\xi',\txi')&=& \frac{\beta_j}{2 \pi} e^{-i
\kappa \xi'}\delta(\txi-\txi'),
\label{eq:A.13}
\end{eqnarray}
which can be solved with the method of characteristics.

We obtain that 
\begin{equation}
  V_{jj}(\om,\kappa,\txi,Z;\xi',\txi') = \frac{\beta_j}{2 \pi} e^{-i
    \kappa \xi'} \delta \left( \txi - \txi'- \kappa Z\right) \exp
    \left[ - \frac{2}{\cS_j}\int_0^Z d s \, {C}_o\Big( \frac{\txi' +
    \kappa s}{\sqrt{\beta_j}}\Big) \right],
\label{eq:A.14}
\end{equation}
and tracing back out transformations (\ref{eq:A.1}), (\ref{eq:A.5}),
(\ref{eq:A.7}) and (\ref{eq:A.12}), we get
\begin{eqnarray}
  \EE\left[ \cT_j(\om,X_1,X_1',Z)
  \overline{\cT_j(\om,X_2,X_2',Z)}\right] = \frac{\beta_j}{2 \pi Z}
  \exp \left[\frac{i \beta_j[(X_1-X_1')^2- (X_2-X_2')^2]}{2 Z}-
  \right. \nonumber \\ \left.  \frac{2}{\cS_j}\int_0^Z d s \, {C}_o\Big[
  (X_1-X_2) \frac{s}{Z} + (X_1'-X_2')\left(1-\frac{s}{Z}\right) \Big]
  \right].
\label{eq:A.15}
\end{eqnarray}
This is the result stated in Proposition \ref{prop.2}.

\subsubsection{Two mode moments}
It is not possible to obtain a closed form solution of (\ref{eq:A.8}),
unless we make further assumptions. We consider the low-SNR regime 
described in section \ref{sect:SNR}, and suppose  that 
\begin{equation}
  |X_1-X_2| \lesssim X_{d,j}(\om) \ll \ell.
\label{eq:AS2M}
\end{equation}
This is the condition under which the diagonal moments $\EE\left[
  \cT_j \overline{\cT_j} \right]$ are not exponentially small, by
Proposition \ref{prop.2ap}. The two mode moments cannot be larger than
the diagonal ones, so they are essentially zero when (\ref{eq:AS2M}) 
does not hold.

Note that in (\ref{eq:A.8}) $\tk$ is the dual variable to $ \txi \sim
\sqrt{\beta_j}(X_1-X_2), $ and that $q$ is in the support of $\widehat
{C}_o$, so $ |q| \le 1/\ell.  $ Therefore,
\[
|\tk| \sim \frac{1}{\sqrt{\beta_j} |X_1-X_2|} \gg \frac{1}{
  \sqrt{\beta_j} \ell} \gtrsim \frac{|q|}{\sqrt{\beta_j}},
\]
and we can expand the Wigner transform in (\ref{eq:A.8}) around $\tk$.
The exponential can also be expanded when
\begin{equation}
\label{eq:AS3M}
\frac{q \xi }{\sqrt{\beta_j}}\frac{2|\sqrt{\beta_j}-\sqrt{\beta_l}|}{
  (\sqrt{\beta_j}+\sqrt{\beta_l})} \lesssim \frac{X}{\ell}
\frac{2|\sqrt{\beta_j}-\sqrt{\beta_l}|}{
  (\sqrt{\beta_j}+\sqrt{\beta_l})} \ll 1,
\end{equation}
meaning that $j$ and $l$ are close. We return at the end of this
section to this point.

Assumptions (\ref{eq:AS2M})-(\ref{eq:AS3M}) justify the approximation
of the right hand side in (\ref{eq:A.8}) by the second-order expansion
in $q$ of the product of the exponential and the Wigner transform. We
obtain that 
\begin{eqnarray}
  \left(\partial_Z + \tk \partial_\xi\right) W_{jl} \approx 
  -\frac{1}{\beta_j \ell^2\sqrt{\cS_j \cS_l}} \left( \frac{\sqrt{\beta_j}+
\sqrt{\beta_l}}{2 \sqrt{\beta_l}}\right)^2 \left[ i \partial_{\tk} - 
\xi \frac{2 (\sqrt{\beta_j}-\sqrt{\beta_l})}{(\sqrt{\beta_j} + 
\sqrt{\beta_l})}\right]^2 W_{jl}
\end{eqnarray}
for $Z>0$, with initial condition (\ref{eq:A.9}). This equation is
solved in \cite{Fannj}. The result follows from the inverse Fourier
transform in $\tk$ of the solution, and from (\ref{eq:A.1}),
(\ref{eq:A.5}),
\begin{eqnarray}
  \EE\left[ \cT_j(\om,X_1,X_1',Z) \overline{\cT_l(\om,X_2,X_2',Z)}\right] 
  \approx \frac{\sqrt{\beta_j \beta_l}}{2 \pi Z} 
  \mbox{\rm  sinc}^{-\frac{1}{2}} 
  \left\{ \frac{(1+i)Z}{\ell} \left[ \frac{\beta_j-\beta_l}{\beta_j \beta_l 
        \sqrt{\cS_j \cS_l}}\right]^{\frac{1}{2}}\right\} \times 
  \nonumber \\ 
  \exp \left\{-\left(
      \frac{1}{\sqrt{\cS_j}}-\frac{1}{\sqrt{\cS_l}}\right)^2 Z + 
\frac{i |(X_1-X_1') \beta_j - (X_2-X_2') \beta_l|^2}{2 Z (\beta_j-\beta_l)} + 
  \right. \nonumber \\ 
  \frac{\beta_j \beta_l \left(|X_1-X_2|^2+|X_1'-X_2'|^2\right)}{
    \ell (\beta_j-\beta_l)(1+i)} 
  \left[ \frac{\beta_j-\beta_l}{\beta_j \beta_l 
      \sqrt{\cS_j \cS_l}}\right]^{\frac{1}{2}} \mbox{cot} \left[ 
    \frac{(1+i)Z}{\ell} \left[ \frac{\beta_j-\beta_l}{\beta_j \beta_l 
        \sqrt{\cS_j \cS_l}}\right]^{\frac{1}{2}}\right] - 
  \nonumber \\
  \left. \frac{2 \beta_j \beta_l (X_1-X_2)(X_1'-X_2')}{\ell 
      (\beta_j-\beta_l)(1+i)} 
    \left[ \frac{\beta_j-\beta_l}{\beta_j \beta_l 
        \sqrt{\cS_j \cS_l}}\right]^{\frac{1}{2}} \sin^{-1} \left[ 
      \frac{(1+i)Z}{\ell} \left[ \frac{\beta_j-\beta_l}{\beta_j \beta_l 
          \sqrt{\cS_j \cS_l}}\right]^{\frac{1}{2}}\right]\right\}
\label{eq:A.16}
\end{eqnarray}
Formula (\ref{eq:A.16}) is complicated, but it can be simplified under 
the assumption that 
\begin{equation}
\frac{Z^2 |\beta_j-\beta_l|}{\beta_j \beta_l \ell^2 \sqrt{\cS_j \cS_l}} 
\ll 1.
\label{eq:A.17}
\end{equation}
Then, we can expand the sinc, cot and $\sin^{-1}$ functions in
(\ref{eq:A.16}) and obtain the simpler formula
\begin{eqnarray}
  \EE\left[ \cT_j(\om,X_1,X_1',Z) \overline{\cT_l(\om,X_2,X_2',Z)}\right] 
  \approx \frac{\sqrt{\beta_j \beta_l}}{2 \pi Z}  \exp \left[
    \frac{i (\beta_j(X_1-X_1')^2 - \beta_l (X_2-X_2')^2)}{2 Z} 
    \right] \times \nonumber \\ 
  \exp \left[ -\Big(
    \frac{1}{\sqrt{\cS_j}}-\frac{1}{\sqrt{\cS_l}}\Big)^2 
    \hspace{-0.05in} Z -  \frac{(X_1-X_2)^2 + (X_1'-X_2')^2 + 
   (X_1-X_2)(X_1'-X_2')}{2 \sqrt{X_{d,j} X_{d,l}}}\right]. 
\label{eq:A.18}
\end{eqnarray}

It remains to justify assumptions (\ref{eq:AS3M}) and (\ref{eq:A.17}). 
Because of the exponential decay in $Z$, we note that the moments are 
essentially zero unless 
\[
\Big(
      \frac{1}{\sqrt{\cS_j}}-\frac{1}{\sqrt{\cS_l}}\Big)^2 
\hspace{-0.05in} Z \lesssim 1.
\]
But in our low-SNR regime this translates to 
\[
\Big(1-\sqrt{\frac{\cS_j}{\cS_l}}\Big)^2 \lesssim
\frac{\cS_j}{\gamma \cS_{1}} \ll 1,
\]
by definition (\ref{eq:gamma}), and it is satisfied only when $j = l$.
This justifies the assumptions, and it means that the modes are
essentially decorrelated.

\subsection{Two frequency moments}
\label{ap:2fM}
The calculation of the two frequency moments is exactly as in the
previous section, with $\beta_j$ replaced by $\beta_j(\om_1)$ and
$\beta_l$ replaced by $\beta_j(\om_2)$. We only consider the case
$j=l$, because the modes decorelate as explained above.  The moment
formula follows from (\ref{eq:A.16}), with $\beta_j$ replaced by
$\beta_j(\om_1)$ and $\beta_l$ replaced by $\beta_j(\om_2)$, and
similar for $\cS_j$ and $\cS_l$.  We can simplify it under the
assumption that $|\om_1 - \om_2|$ is sufficiently small to make
 first-order expansions in $\om_1-\om_2$.  Let $\tom$ and $\om$ be the center
and difference frequencies
\[
\tom = \om_1-\om_2, \qquad \om = \frac{\om_1+\om_2}{2}.
\]
We have from (\ref{eq:A.16}) that 
\begin{eqnarray}
  \EE\left[ \cT_j\left(\om+\frac{ \tom}{2},X_1,X_1',Z\right)
\overline{\cT_l\left(\om-\frac{ \tom}{2},X_2,X_2',Z\right)}\right]
\approx \frac{\beta_j(\om)}{2 \pi Z} \mbox{\rm sinc}^{-\frac{1}{2}}
\left\{ \frac{(1+i)Z}{\ell\beta_j(\om)} \left[ \frac{\tom
\partial_{\om} \beta_j(\om)}{\cS_j(\om)}\right]^{\frac{1}{2}} \right\}
\times \nonumber \\ \exp \left\{-\tom^2 \Big(\partial_\om
\frac{1}{\sqrt{\cS_j(\om)}}\Big)^2
    \hspace{-0.05in} Z + + \frac{i |[(X_1-X_1')-(X_2-X_2')]
    \beta_j(\om) + \tom \frac{[(X_1-X_1')+ (X_2-X_2')]}{2}
    \beta'_j(\om)|^2}{2 Z \tom \beta'_j(\om)}
    + \right. \nonumber \\ \frac{\beta_j(\om)[|X_1-X_2|^2 +
    |X_1'-X_2'|^2]}{\ell \tom \beta'_j(\om) (1+i)}\left[
    \frac{\tom \beta'_j(\om)}{
    \cS_j(\om)}\right]^{\frac{1}{2}} \mbox{cot} \left[
    \frac{(1+i)Z}{\beta_j(\om)\ell} \left[ \frac{\tom 
    \beta'_j(\om)}{ \cS_j(\om)}\right]^{\frac{1}{2}}\right] - \nonumber
    \\ \left. \frac{2 \beta_j(\om) (X_1-X_2)(X_1'-X_2')}{\ell \tom
    \beta'_j(\om) (1+i)}\left[ \frac{\tom 
    \beta'_j(\om)}{ \cS_j(\om)}\right]^{\frac{1}{2}} \sin^{-1} \left[
    \frac{(1+i)Z}{\beta_j(\om)\ell} \left[ \frac{\tom 
    \beta'_j(\om)}{ \cS_j(\om)}\right]^{\frac{1}{2}}\right]\right\}.
    \quad
\label{eq:A.20}
\end{eqnarray}

\subsection{Frequency decorrelation}
\label{ap:freqDec}
To study the decorrelation over frequency offsets, let $X_1 = X_2$ and
$X_1' = X_2'$ in (\ref{eq:A.20})
\begin{eqnarray}
  \EE\left[ \cT_j\Big(\om + \frac{\tom}{2},X,X',Z\Big) \overline{\cT_j
      \Big(\om-\frac{\tom}{2},X,X',Z\Big)}\right] \approx
      \frac{\beta_j(\om)}{2 \pi Z} \mbox{\rm sinc}^{-\frac{1}{2}}
      \left\{ \frac{(1+i)Z}{\ell\beta_j(\om)} \left[ \frac{\tom
      \partial_{\om} \beta_j(\om)}{\cS_j(\om)}\right]^{\frac{1}{2}}
      \right\} \nonumber \\ \times \exp \left\{-\left[\tom
      \partial_\om \cS_j^{-\frac{1}{2}}(\om)\right]^2
    \hspace{-0.05in} Z +  \frac{i \tom \beta'_j(\om)
      (X-X')^2}{2Z}\right\}.
\label{eq:2M.7}
\end{eqnarray}
We have two factors that decay exponentially in $\tom$. The first is
the sinc, decaying at the rate 
\begin{equation}
  |\tom| \ll \Omega_{d,j}(\om) = \frac{\cS_j(\om) \beta_j^2(\om)
    \ell^2}{Z^2 |\beta'_j(\om)|} =
    \frac{\beta_j(\om)}{|\beta'_j(\om)|}\frac{\cS_j(\om)
    \beta_j(\om)\ell^2}{\gamma^2 \cS_1^2(\om)}.
\label{eq:2M.3a}
\end{equation}
and the second is the Gaussian with standard deviation
\begin{equation}
  \Omega_j(\om) = \frac{1}{\sqrt{2 Z} \left|\partial_\om
      \cS_j^{-\frac{1}{2}}(\om)\right|} = \frac{\beta_j(\om)}{
    |\beta'_j(\om)|} \sqrt{ \frac{\cS_j(\om)}{2
\gamma \cS_1(\om)}}.
\label{eq:2M.8a}
\end{equation}
Note that
\begin{equation}
  \frac{|\beta'_j(\om)|}{\beta_j(\om)} =
  \frac{\left(N+\alpha(\om)-{1}/{2}\right)^2}{\om \left[
  \left(N+\alpha(\om)-{1}/{2}\right)^2 -
  \left(j-{1}/{2}\right)^2\right]},
\end{equation}
and using equations (\ref{eq:HOM.13}), (\ref{eq:DefAlpha}),
(\ref{eq:STAT.3}), and the high-frequency assumption $N \gg 1$, we
have
\begin{equation}
  \Omega_{d,j}(\om) \approx \frac{\om \ell^2 \beta_1(\om)} {16
  \gamma^2 \cS_1(\om)} \frac{\left[
  \left(N+\alpha(\om)-\frac{1}{2}\right)^2 -
  \left(j-\frac{1}{2}\right)^2\right]^{5/2}}{N^5 (j-1/2)^4}.
\end{equation}
Here
\begin{equation}
  \frac{\ell^2 \beta_1}{\cS_1(\om)} \approx \frac{\sigma^2 }{32 N}
  \left(\frac{\pi \ell }{{{\mathcal D}}}\right)^3
  \approx \frac{\sigma^2 (\ell k)^3}{32 N^4} =
  \frac{\pi^3 \sigma^2 (\ell/\lambda)^3}{4 N^4},
\end{equation}
and $\ell/\lambda = O(1)$, because the scaled correlation length is
similar to the wavelength $\lambda$.  Moreover, the rate
(\ref{eq:2M.8a}) is given by
\begin{equation}
  \Omega_j(\om) \approx \frac{\om}{4 \sqrt{2\gamma}} \frac{\left[
      \left(N+\alpha(\om)-\frac{1}{2}\right)^2 -
      \left(j-\frac{1}{2}\right)^2\right]^{3/2}}{N^3 (j-1/2)^2},
\end{equation}
and it is larger than $\Omega_{d,j}(\om)$.

Thus, we call $\Omega_{d,j}(\om)$ the mode-dependent {\em decoherence
  frequency}, the frequency scale over which the second moments decay.
Note that when the frequency offsets satisfy $|\tom| \ll \Omega_{d,j}$
the moment formula (\ref{eq:2M.7}) simplifies to  expression
(\ref{eq:2M.4}) in Proposition \ref{prop.3}, because
\[
\exp \left[-\frac{\tom^2}{2 \Omega_j^2(\om)}\right] \approx 1,
\] 
when
\[ |\tom| \ll \Omega_{d,j} \ll \Omega_j(\om). \]

\section{The fourth moments}
\label{ap:FM}
We denote the  moments by 
\begin{equation}
\label{eq:FM.1}
M_{jlJL}:=\EE \left[ \cT_j(\om_1,X_1,X_1',Z) \overline{\cT_l(\om_2,X_2,X_2',Z)}
  \cT_J(\om_3,Y_1,Y_1',Z) \overline{\cT_L(\om_4,Y_2,Y_2',Z)}\right],
\end{equation}
and obtain from (\ref{eq:LR.14}) that they 
satisfy the partial differential equation
\begin{eqnarray}
\label{eq:FM.2}
\partial_Z M_{jlJL} &=& \left[\frac{i}{2 \beta_j} \partial_{X_1}^2 
  -\frac{i}{2 \beta_l} \partial_{X_2}^2 +\frac{i}{2 \beta_J} \partial_{Y_1}^2 
  -\frac{i}{2 \beta_L} \partial_{Y_2}^2 \right]  M_{jlJL}  \nonumber \\
&&+ \left[- \Big(\frac{1}{\sqrt{\cS_j}}-
  \frac{1}{\sqrt{\cS_l}}\Big)^2 - \Big(\frac{1}{\sqrt{\cS_J}}-
  \frac{1}{\sqrt{\cS_L}}\Big)^2 - \frac{2 {C}_o(X_1-X_2)}{\sqrt{\cS_j \cS_l}} - 
  \frac{2 {C}_o(Y_1-Y_2)}{\sqrt{\cS_J \cS_L}}\right]M_{jlJL} \nonumber \\
&& + \frac{2}{R_o(0)} \left[
  \frac{R_o(X_1-Y_2)}{\sqrt{\cS_j \cS_L}}-\frac{R_o(X_1-Y_1)}{\sqrt{\cS_j \cS_J}}
  -\frac{R_o(X_2-Y_2)}{\sqrt{\cS_l \cS_L}}
  +\frac{R_o(X_2-Y_1)}{\sqrt{\cS_l \cS_J}}\right] M_{jlJL},
\label{eq:FM.3}
\end{eqnarray}
for $Z>0$, with the initial condition 
\begin{equation}
  M_{jlJL} = \delta(X_1-X_1') \delta(X_2-X_2') \delta(Y_1-Y_1') 
\delta(Y_2-Y_2'), \quad \mbox{at} 
  ~ Z = 0.
\label{eq:FM.4}
\end{equation}

Let us consider the case $j=l$, $J = L$, $\om_1 = \om_2 = \om$, and
$\om_3 = \om_4 = \om'$. These moments $M_{jjJJ}$ are needed in
section \ref{sect:TR} to show the statistical stability of the time
reversal function in the case of an array that spans the entire depth
of the waveguide. 
We look for the fourth-order moment for $X_1=X_2$ and $Y_1=Y_2$ 
in the support of the array. 
So we parameterize
\begin{eqnarray}
X_1 &=& |A_X| \xi +  X_{d,j}(\omega) \frac{u}{2}, \quad  
X_2 =|A_X|  \xi - X_{d,j}(\omega) \frac{u}{2}, \\
X_1' &=& |A_X|  \xi' + X_{d,j}(\omega) \frac{u'}{2}, \quad  
X_2' = |A_X|  \xi' -  X_{d,j}(\omega)\frac{u'}{2}, \label{eq:FM.5}\\
Y_1 &=& |A_X|  \zeta +  X_{d,J}(\omega') \frac{v}{2}, \quad  
Y_2 = |A_X| \zeta -X_{d,J}(\omega')  \frac{v}{2}, \\
Y_1' &=& |A_X|  \zeta' +  X_{d,J}(\omega') \frac{v'}{2}, \quad  
Y_2' = |A_X|  \zeta' - X_{d,J}(\omega') \frac{v'}{2}. \label{eq:FM.6}
\end{eqnarray}
Equation (\ref{eq:FM.3}) becomes (remember ${C}_o''(0) = 1/\ell^2$)
\begin{eqnarray}
\nonumber \partial_Z M_{jjJJ} \approx \left[ \frac{i}{\beta_j X_{d,j}
|A_X|} \partial_\xi \partial_u + \frac{i}{\beta_J X_{d,J} |A_X|}
\partial_\zeta \partial_v - \frac{X_{d,j}^2 u^2}{\ell^2
\cS_j}-\frac{X_{d,J}^2 v^2}{ \ell^2 \cS_J} \right.\\ \left. - \frac{2
X_{d,j}X_{d,J}uv {C}_o''[|A_X|(\xi-\zeta)]}{ \sqrt{\cS_j\cS_J}}
\right]M_{jjJJ},
\label{eq:FM.8}
\end{eqnarray}
with the initial condition
\begin{equation}
\label{eq:FM.9}
M_{jjJJ} = \frac{1}{X_{d,j} X_{d,j} |A_X|^2} \delta(u-u')
\delta(\xi-\xi') \delta(v-v') \delta(\zeta-\zeta'), \qquad \mbox{at}~
~ Z = 0.
\end{equation}

We address two cases: 

\vspace{0.05in} \noindent \textbf{Case 1:} The array diameter $|A_X|$
is much larger than $\ell$. This allows us to simplify Equation
(\ref{eq:FM.8}) as
\begin{eqnarray}
\partial_Z M_{jjJJ} &\approx& \left[ \frac{i}{\beta_j X_{d,j} |A_X|}
\partial_\xi \partial_u + \frac{i}{\beta_J X_{d,J} |A_X|}
\partial_\zeta \partial_v - \frac{X_{d,j}^2 u^2}{\ell^2
\cS_j}-\frac{X_{d,J}^2 v^2}{ \ell^2 \cS_J} \right]M_{jjJJ},
\label{eq:FM.8b}
\end{eqnarray}
which has a separable form in $(u,\xi)$ and $(v,\zeta)$, and we get
(following the same method as in the case of second-order moments):
\begin{eqnarray}
\nonumber M_{jjJJ}& \approx& \frac{\beta_j \beta_J}{4 \pi^2 Z^2} \exp
\left[ - \frac{i \beta_jX_{d,j} |A_X| (\xi-\xi') (u-u')}{Z} - \frac{i
\beta_JX_{d,J} |A_X| (\zeta-\zeta')( v-v')}{Z}\right]\\ &&\times \exp
\left\{ - \frac{Z}{3 \ell^2} \left[ \frac{X_{d,j}^2
(u^2+u'^2+uu')}{\cS_j} + \frac{X_{d,J}^2 (v^2+v'^2+vv')}{\cS_J}
\right] \right\} .
\end{eqnarray}
Equivalently, in terms of the original variables, 
\begin{eqnarray}
\nonumber M_{jjJJ} &\approx&\frac{\beta_j \beta_J}{4 \pi^2 Z^2} \exp
\left[ - \frac{i \beta_j [ (X_1-X_1')^2-(X_2-X_2')^2] }{2Z} - \frac{i
\beta_J [(Y_1-Y_1')^2-(Y_2-Y_2')^2]}{2 Z} \right]\\ \nonumber &&\times
\exp \left[ - \frac{1}{2} \frac{(X_1-X_2)^2+(X_1'-X_2')^2
+(X_1-X_2)(X_1'-X_2')}{X_{d,j}^2} \right]\\ &&\times \exp\left[ -
\frac{1}{2}
\frac{(Y_1-Y_2)^2+(Y_1'-Y_2')^2+(Y_1-Y_2)(Y_1'-Y_2')}{X_{d,J}^2}
\right],
\end{eqnarray}
which is equal to 
$  \EE\left[\cT_j \overline{\cT_j}\right]\EE\left[\cT_J \overline{\cT_J}
  \right] $.
  
\vspace{0.05in} \noindent \textbf{Case 2:} The array diameter $|A_X|$
is smaller than $\ell$. Then Equation (\ref{eq:FM.3}) becomes
\begin{eqnarray}
\partial_Z M_{jjJJ} &\approx& \left[ \frac{i}{\beta_j X_{d,j} |A_X|}
\partial_\xi \partial_u + \frac{i}{\beta_J X_{d,J} |A_X|}
\partial_\zeta \partial_v - \frac{X_{d,j}^2 u^2}{\ell^2
\cS_j}-\frac{X_{d,J}^2 v^2}{ \ell^2 \cS_J}- \frac{2 X_{d,j}X_{d,J}uv
}{\ell \sqrt{\cS_j\cS_J}} \right]M_{jjJJ}, \qquad 
\label{eq:FM.8c}
\end{eqnarray}
with the initial condition (\ref{eq:FM.9}).  This equation can be
solved explicitly after Fourier transforming in $\xi$ and $\zeta$.  If
we let
\begin{equation}
\hat M_{jjJJ} = \int d \xi \int d \zeta M_{jjJJ} e^{ i K_\xi \xi
+ i K_\zeta \zeta},
\end{equation}
then we have 
\begin{equation}
\left[ \partial_Z + \frac{K_\xi}{\beta_j X_{d,j} |A_X|} \partial_u +
\frac{K_\zeta}{\beta_J X_{d,J} |A_X| } \partial_v\right] \hat M_{jjJJ}
\approx -\Big(\frac{X_{d,j} u}{\ell \sqrt{\cS_j}}+\frac{X_{d,J} v}{
\ell \sqrt{\cS_J}}\Big)^2 \hat M_{jjJJ},
\end{equation}
for $Z >0$, and 
\begin{equation}
\hat M_{jjJJ} = \frac{1}{X_{d,j} X_{d,j} |A_X|^2}\delta(u-u')
 \delta(v-v') e^{i K_\xi \xi'+i K_\zeta \zeta'} , \qquad \mbox{at}~ ~
 Z = 0.
\end{equation}
The solution is given by the method of characteristics 
\begin{eqnarray}
\nonumber \hat M_{jjJJ} \approx \delta\left(u- u'- \frac{K_\xi
  Z}{\beta_j X_{d,j} |A_X|}\right)\delta\left(v-v'- \frac{K_\zeta
  Z}{\beta_J X_{d,J} |A_X|}\right) \exp \left[ -\frac{Z^3}{3 \ell^2}
  \Big(\frac{K_\xi}{\beta_j \sqrt{\cS_j}|A_X|} +\frac{K_\zeta}{\beta_J
  \sqrt{\cS_J}|A_X|}\Big)^2 \right. \\ \left. -\frac{Z^2}{\ell^2}
  \Big(\frac{K_\xi}{\beta_j \sqrt{\cS_j} |A_X|}
  +\frac{K_\zeta}{\beta_J \sqrt{\cS_J}|A_X|}\Big)
  \Big(\frac{X_{d,j}u'}{\sqrt{\cS_j}}
  +\frac{X_{d,J}v'}{\sqrt{\cS_J}}\Big) -\frac{Z}{\ell^2}
  \Big(\frac{X_{d,j}u'}{\sqrt{\cS_j}}
  +\frac{X_{d,J}v'}{\sqrt{\cS_J}}\Big)^2 \right], \quad \quad \quad
\end{eqnarray}
and the moment estimate follows from the inverse Fourier transform,
\begin{eqnarray*}
\nonumber &&M_{jjJJ} \approx \frac{\beta_j \beta_J}{4 \pi^2 Z^2} \exp
\left[ - \frac{i \beta_j |A_X|X_{d,j} (\xi-\xi') (u-u')}{Z} - \frac{i
\beta_J |A_X|X_{d,J} (\zeta-\zeta')( v-v')}{Z}\right]\\ &&\times \exp
\left\{ - \frac{Z}{3 \ell^2} \left[ \Big( \frac{u
X_{d,j}}{\sqrt{\cS_j}} + \frac{v X_{d,J}}{\sqrt{\cS_J}}\Big)^2 +\Big(
\frac{u'X_{d,j}}{\sqrt{\cS_j}} + \frac{v'
X_{d,J}}{\sqrt{\cS_J}}\Big)^2 +\Big( \frac{uX_{d,j}}{\sqrt{\cS_j}} +
\frac{v X_{d,J}}{\sqrt{\cS_J}}\Big)\Big(
\frac{u'X_{d,j}}{\sqrt{\cS_j}} + \frac{v'X_{d,J}}{\sqrt{\cS_J}}\Big)
\right] \right\} .
\end{eqnarray*}
Equivalently, in terms of the original variables, 
\begin{eqnarray}
\nonumber M_{jjJJ} &\approx& \frac{\beta_j \beta_J}{4 \pi^2 Z^2} \exp
\left[ - \frac{i \beta_j [ (X_1-X_1')^2-(X_2-X_2')^2] }{2Z} - \frac{i
\beta_J [(Y_1-Y_1')^2-(Y_2-Y_2')^2]}{2 Z} \right]\\ \nonumber &&\times
\exp \left[ - \frac{1}{2}\Big( \frac{X_1-X_2}{X_{d,j}} +
\frac{Y_1-Y_2}{X_{d,J}}\Big)^2 - \frac{1}{2}\Big(
\frac{X_1'-X_2'}{X_{d,j}} + \frac{Y_1'-Y_2'}{X_{d,J}}\Big)^2
\right. \\ && \quad\quad\quad \left.  - \frac{1}{2}\Big(
\frac{X_1-X_2}{X_{d,j}} + \frac{Y_1-Y_2}{X_{d,J}}\Big)\Big(
\frac{X_1'-X_2'}{X_{d,j}} + \frac{Y_1'-Y_2'}{X_{d,J}}\Big) \right].
\end{eqnarray}
If $X_1=X_2$ and $Y_1=Y_2$, then
\begin{eqnarray}
\nonumber M_{jjJJ} &\approx& \frac{\beta_j \beta_J}{4 \pi^2 Z^2} \exp
\left[ - \frac{i \beta_j [ (X_1-X_1')^2-(X_2-X_2')^2] }{2Z} - \frac{i
\beta_J [(Y_1-Y_1')^2-(Y_2-Y_2')^2]}{2 Z} \right]\\ \nonumber &&\times
\exp \left[ - \frac{1}{2} \Big( \frac{X_1'-X_2'}{X_{d,j}}+
\frac{Y_1'-Y_2'}{X_{d,J}} \Big)^2 \right],
\end{eqnarray}
while
\begin{eqnarray}
\nonumber \EE\left[\cT_j \overline{\cT_j}\right]\EE\left[\cT_J
  \overline{\cT_J} \right] &\approx& \frac{\beta_j \beta_J}{4 \pi^2
  Z^2} \exp \left[ - \frac{i \beta_j [ (X_1-X_1')^2-(X_2-X_2')^2]
  }{2Z} - \frac{i \beta_J [(Y_1-Y_1')^2-(Y_2-Y_2')^2]}{2 Z} \right]\\
  \nonumber &&\times \exp \left[ - \frac{1}{2} \left(
  \frac{(X_1'-X_2')^2}{X_{d,j}^2}+ \frac{(Y_1'-Y_2')^2}{X_{d,J}^2}
  \right) \right].
\end{eqnarray}
Here we can see that the fourth-order moment is not equal to the
product of the second-order moments.

 \end{document}